\documentclass[11pt]{article} 
\usepackage{color}
\usepackage{graphicx}
\usepackage{amsmath}
\usepackage{amssymb}
\usepackage{theorem}
\usepackage{euscript}
\usepackage[usenames,dvipsnames]{pstricks}
\usepackage{pstricks-add}
\usepackage{pst-grad}
\usepackage{pst-plot} 
\usepackage{pst-eucl} 
\usepackage{theorem}
\newtheorem{plctheorem}{Theorem}[section]

\newtheorem{proposition}[plctheorem]{Proposition}
\newtheorem{definition}[plctheorem]{Definition}

\theoremstyle{plain}{\theorembodyfont{\rmfamily}
}
\theoremstyle{plain}{\theorembodyfont{\rmfamily}
}
\theoremstyle{plain}{\theorembodyfont{\rmfamily}
\newtheorem{algorithm}[plctheorem]{Algorithm}}
\theoremstyle{plain}{\theorembodyfont{\rmfamily}
\newtheorem{plcexample}[plctheorem]{Example}}

\theoremstyle{plain}{\theorembodyfont{\rmfamily}
}
\theoremstyle{plain}{\theorembodyfont{\rmfamily}
\newtheorem{problem}[plctheorem]{Problem}}

\newcounter{numperceuse}
\renewcommand\thenumperceuse{\roman{numperceuse}}
\newcounter{numscie}
\renewcommand\thenumscie{\roman{numscie}}

\begin{document}
\title{\sffamily Proximal Splitting Methods in Signal 
Processing\thanks{This work was supported by the Agence Nationale 
de la Recherche under grants ANR-08-BLAN-0294-02 and 
ANR-09-EMER-004-03.}}
\author{Patrick L. Combettes\thanks{UPMC Universit\'e Paris 06,
Laboratoire Jacques-Louis Lions -- UMR CNRS 7598, 75005 Paris, 
France} \;\:and 
Jean-Christophe Pesquet\thanks{Laboratoire d'Informatique Gaspard 
Monge, UMR CNRS 8049, Universit\'e Paris-Est, 77454 Marne la Vall\'ee 
Cedex 2, France.}}

\date{~}

\maketitle

\medskip

\abstract{The proximity operator of a convex function is a 
natural extension of the notion of a projection operator onto a convex
set. This tool, which plays a central role in the analysis and the 
numerical solution of convex optimization problems, has recently been
introduced in the arena of inverse problems and, especially, in
signal processing, where it has become 
increasingly important. In this paper, we review the basic
properties of proximity operators which are relevant to signal
processing and present optimization methods based on 
these operators. 
These proximal splitting methods are shown to capture and 
extend several well-known algorithms in a unifying framework. 
Applications of proximal methods in signal recovery and
synthesis are discussed.
}

{\bfseries Keywords.}
Alternating-direction method of multipliers,
backward-backward algorithm,
convex optimization, 
denoising, 
Douglas-Rachford algorithm, 
forward-backward algorithm, 
frame, 
Landweber method, 
iterative thresholding, 
parallel computing, 
Peaceman-Rachford algorithm, 
proximal algorithm, 
restoration and reconstruction, 
sparsity,
splitting. 

\bigskip
\noindent{\bfseries AMS 2010 Subject Classification:} 
90C25, 65K05, 90C90, 94A08

\section{Introduction}
\label{plcsec:1}

Early signal processing methods were essentially linear, as they 
were based on classical functional analysis and linear algebra. 
With the development of nonlinear analysis in mathematics in 
the late 1950s and early 1960s (see the bibliographies
of \cite{plcAubi98,plcZeidXX}) and the 
availability of faster computers, nonlinear techniques have 
slowly become prevalent. In particular, convex optimization 
has been shown to provide efficient 
algorithms for computing reliable solutions in a
broadening spectrum of applications.

Many signal processing problems can \emph{in fine} be 
formulated as convex optimization problems of the form 
\begin{equation}
\label{plce:1}
\underset{x\in{\mathbb R}^N}{\mathrm{minimize}}\;\;f_1(x)+\cdots+f_m(x),
\end{equation}
where $f_1,\ldots,f_m$ are convex functions from ${\mathbb R}^N$ to 
$\,\left]-\infty,+\infty\right]$.
A major difficulty that arises in solving this problem stems from the
fact that, typically, some of the functions are not differentiable,
which rules out conventional smooth optimization techniques.
In this paper, we describe a class of efficient convex optimization
algorithms to solve \eqref{plce:1}. These methods proceed by
\emph{splitting} in that the functions $f_1$, \ldots, $f_m$ are 
used individually so as to yield an easily implementable algorithm. 
They are called \emph{proximal} because
each nonsmooth function in \eqref{plce:1} is involved via its proximity 
operator. Although proximal methods, which can be traced back to the 
work of Martinet \cite{plcMart70}, have been introduced in signal
processing only recently \cite{plcSmai01,plcSmms05}, their use 
is spreading rapidly.

Our main objective is to familiarize the 
reader with proximity operators, their main properties, and a variety
of proximal algorithms for solving signal and image 
processing problems. The power and flexibility of proximal methods 
will be emphasized. 
In particular, it will be shown that a number of apparently
unrelated, well-known algorithms (e.g., iterative thresholding,
projected Landweber, projected gradient, alternating projections, 
alternating-direction method of multipliers, 
alternating split Bregman) are special instances of proximal 
algorithms. In this respect, the proximal formalism provides a 
unifying framework for analyzing and developing a broad class of 
convex optimization algorithms. 
Although many of the subsequent results are extendible to 
infinite-dimensional spaces, we restrict ourselves to a 
finite-dimensional setting to avoid technical digressions.

The paper is organized as follows.
Proximity operators are introduced in Section~\ref{plcsec:2}, where we
also discuss their main properties and provide examples.
In Sections~\ref{plcsec:3} and~\ref{plcsec:4}, 
we describe the main proximal
splitting algorithms, namely the forward-backward algorithm and the
Douglas-Rachford algorithm. In Section~\ref{plcsec:5}, we present a
proximal extension of Dykstra's projection method which is 
tailored to problems featuring strongly convex objectives.
Composite problems involving linear transformations of the variables
are addressed in Section~\ref{plcsec:6}. The algorithms discussed so far
are designed for $m=2$ functions. 
In Section~\ref{plcsec:7}, we discuss parallel variants 
of these algorithms
for problems involving $m\geq 2$ functions.
Concluding remarks are given in Section~\ref{plcsec:9}. 

{\bfseries Notation.}
We denote by ${\mathbb R}^N$ the usual $N$-dimensional Euclidean
space, by $\|\cdot\|$ its norm, and by $I$ the identity matrix. 
Standard definitions and notation from 
convex analysis will be used \cite{plcLivr10,plcHiri93,plcRock70}.
The domain of a function 
$f\colon{\mathbb R}^N\to\,\left]-\infty,+\infty\right]$ is
${\mathrm{dom}}\, f=\{x\in{\mathbb R}^N\:|\:f(x)<+\infty\}$.
$\Gamma_0({\mathbb R}^N)$ is the class of lower 
semicontinuous convex functions 
from ${\mathbb R}^N$ to $\,\left]-\infty,+\infty\right]$ such that 
${\mathrm{dom}}\, f\neq\varnothing$. 
Let $f\in\Gamma_0({\mathbb R}^N)$. The conjugate of $f$ is
the function $f^*\in\Gamma_0({\mathbb R}^N)$ defined by
\begin{equation}
f^*\colon{\mathbb R}^N\to\,\left]-\infty,+\infty\right]\colon u\mapsto
\sup_{x\in{\mathbb R}^N}{x}^\top\!{u}-f(x),
\end{equation}
and the subdifferential of $f$ is the set-valued operator
\begin{equation}
\label{e:subdiff}
\partial f\colon{\mathbb R}^N\to 2^{{\mathbb R}^N}\colon x\mapsto
\big\{u\in{\mathbb R}^N\:|\:(\forall y\in{\mathbb R}^N)\;\:
{(y-x)}^\top{u}+f(x)\leq f(y)\big\}.
\end{equation}
Let $C$ be a nonempty subset of ${\mathbb R}^N$. The indicator 
function of $C$ is
\begin{equation}
\label{e:iota}
\iota_C\colon x\mapsto
\begin{cases}
0,&\text{if}\;\;x\in C;\\
+\infty,&\text{if}\;\;x\notin C,
\end{cases}
\end{equation}
the support function of $C$ is 
\begin{equation}
\label{e:support}
\sigma_C=\iota_C^*\colon{\mathbb R}^N\to\,\left]-\infty,+\infty\right]
\colon u\mapsto\underset{x\in C}{\text{\rm sup\,}}{u}^\top\!{x},
\end{equation}
the distance from $x\in{\mathbb R}^N$ to $C$ is 
$d_C(x)={\text{\rm inf}}_{y\in C}\|x-y\|$, 
and the relative interior of $C$ (i.e., interior of $C$
relative to its affine hull) is the nonempty set
denoted by ${\mathrm{ri}}\, C$.
If $C$ is closed and convex, the projection of $x\in {\mathbb R}^N$ onto 
$C$ is the unique point $P_C x\in C$ such that $d_C(x)=\|x-P_Cx\|$.

\section{From projection to proximity operators}
\label{plcsec:2}
One of the first widely used convex optimization 
splitting algorithms in signal processing is POCS 
(Projection Onto Convex Sets) \cite{plcCens97,plcProc93,plcYoul82}.
This algorithm is employed to recover/synthesize a signal 
satisfying simultaneously several convex constraints.
Such a problem can be formalized within the framework of \eqref{plce:1}
by letting each function $f_i$ be the indicator function of 
a nonempty closed convex set $C_i$ modeling a constraint.
This reduces \eqref{plce:1} to the classical \emph{convex feasibility
problem} 
\cite{plcCens97,plcProc93,plcAiep96,plcHerm80,plcLent81,%
plcStar87,plcStar98,plcTrus84,plcYoul82}
\begin{equation}
\label{plce:feas}
\text{find}\;\;x\in\bigcap_{i=1}^m C_i.
\end{equation}
The POCS algorithm \cite{plcBreg65,plcYoul82} activates each set $C_i$
individually by means of its projection operator $P_{C_i}$. It
is governed by the updating rule
\begin{equation}
\label{plce:pocs}
x_{n+1}=P_{C_1}\cdots P_{C_m}x_n.
\end{equation}
When $\bigcap_{i=1}^m C_i\neq\varnothing$ 
the sequence $(x_n)_{n\in{\mathbb N}}$ thus
produced converges to a solution to \eqref{plce:feas} \cite{plcBreg65}.
Projection algorithms have been enriched with many 
extensions of this basic iteration to solve \eqref{plce:feas}
\cite{plcBaus96,plcSign94,plcImag97,plcLopu97}. 
Variants have also been proposed
to solve more general problems, e.g., that of finding the projection 
of a signal onto an intersection of convex sets 
\cite{plcBoyl86,plcSign03,plcYama98}. 
Beyond such problems, however, projection 
methods are not appropriate and more general operators are required 
to tackle \eqref{plce:1}. Among the various generalizations of the 
notion of a convex projection operator that exist 
\cite{plcBaus96,plcMoor01,plcAiep96,plcLopu97}, proximity operators 
are best suited for our purposes.

The projection $P_{C}x$ of $x\in{\mathbb R}^N$ onto the nonempty closed
convex set $C\subset{\mathbb R}^N$ is the solution to the problem 
\begin{equation}
\label{plce:proj2}
\underset{y\in{\mathbb R}^N}{\mathrm{minimize}}\;\;
\iota_{C}(y)+\frac12\|x-y\|^2.
\end{equation}
Under the above hypotheses, the function $\iota_{C}$ belongs to 
$\Gamma_0({\mathbb R}^N)$. In 1962, Moreau \cite{plcMor62b} 
proposed the following extension of the notion of a projection operator,
whereby the function $\iota_{C}$ in \eqref{plce:proj2} is replaced 
by an arbitrary function $f\in\Gamma_0({\mathbb R}^N)$.

\begin{definition}[Proximity operator]
Let $f\in\Gamma_0({\mathbb R}^N)$. For every $x\in{\mathbb R}^N$, 
the minimization problem
\begin{equation}
\label{plce:prox1}
\underset{y\in{\mathbb R}^N}{\mathrm{minimize}}\;\;f(y)+\frac12\|x-y\|^2
\end{equation}
admits a unique solution, which is denoted by ${\mathrm{prox}}_{f}x$. 
The operator ${\mathrm{prox}}_{f}\colon{\mathbb R}^N\to{\mathbb R}^N$ 
thus defined is the
\emph{proximity operator} of $f$.
\end{definition}

Let $f\in\Gamma_0({\mathbb R}^N)$. The proximity operator of $f$ 
is characterized by the inclusion
\begin{equation}
\label{e:prox1}
(\forall (x,p)\in{\mathbb R}^N\times{\mathbb R}^N)\quad 
p={\mathrm{prox}}_f\,x\quad
\Leftrightarrow\quad x-p\in\partial f(p),
\end{equation}
which reduces to
\begin{equation}
\label{e:prox2}
(\forall (x,p)\in{\mathbb R}^N\times{\mathbb R}^N)\quad p
={\mathrm{prox}}_f\,x\quad
\Leftrightarrow\quad x-p=\nabla f(p)
\end{equation}
if $f$ is differentiable.
Proximity operators have very attractive properties that make 
them particularly well suited for iterative minimization algorithms. 
For instance, ${\mathrm{prox}}_{f}$ is firmly nonexpansive, i.e.,
\begin{multline}
(\forall x\in{\mathbb R}^N)(\forall y\in{\mathbb R}^N)\quad
\|{\mathrm{prox}}_fx-{\mathrm{prox}}_fy\|^2+\|(x-{\mathrm{prox}}_fx)
-(y-{\mathrm{prox}}_fy)\|^2\\
\leq\|x-y\|^2,
\end{multline}
and its fixed point set is precisely the set of minimizers of $f$. 
Such properties allow us to envision the possibility of
developing algorithms based on the proximity operators 
$({\mathrm{prox}}_{f_i})_{1\leq i\leq m}$ to 
solve \eqref{plce:1}, mimicking to 
some extent the way convex feasibility algorithms employ the 
projection operators $(P_{C_i})_{1\leq i\leq m}$ to solve 
\eqref{plce:feas}. As shown in Table~\ref{plct:prop}, 
proximity operators enjoy many additional properties.
One will find in Table~\ref{plct:real} closed-form expressions of
the proximity operators of various functions in $\Gamma_0({\mathbb R})$
(in the case of functions such as $|\cdot|^p$, proximity operators
implicitly appear in several places, e.g., 
\cite{plcAnto01,plcAnto02,plcCham98}).

From a signal processing perspective, proximity 
operators have a very natural interpretation in terms of denoising. 
Let us consider the standard denoising problem of recovering 
a signal $\overline{x}\in{\mathbb R}^N$ from an observation
\begin{equation}
\label{plce:model1}
y=\overline{x}+w,
\end{equation}
where $w\in{\mathbb R}^N$ models noise. This problem
can be formulated as \eqref{plce:prox1}, where $\|\cdot-y\|^2/2$ 
plays the role of a data fidelity term and where $f$ models 
a priori knowledge about $\overline{x}$. Such a formulation derives in 
particular from a Bayesian approach to denoising 
\cite{plcBoum93,plcThom93,plcTitt85} in the presence of Gaussian 
noise and of a prior with a log-concave density $\exp(-f)$.

\refstepcounter{numperceuse}
\begin{table}[htb]
\centering
\caption{Properties of proximity operators
\cite{plcLuis09,plcChau07,plcSiop07,plcInvp08,plcSmms05,plcMore65}:
$\varphi\in\Gamma_0({\mathbb R}^N)$; $C\subset{\mathbb R}^N$ is nonempty, 
closed, and convex; $x\in{\mathbb R}^N$.}
\label{plct:prop}
\begin{tabular}{|p{0.4cm} l|l|l|}
\hline
& Property & $f(x)$ & ${\mathrm{prox}}_{f}\,x$\\
\hline
\hline
\thenumperceuse\refstepcounter{numperceuse} &
translation & $\varphi(x-z)$, $z \in {\mathbb R}^N$ & 
$z+{\mathrm{prox}}_{\varphi}(x-z)$\\
\hline
\thenumperceuse\refstepcounter{numperceuse} &
scaling & $ \varphi(x/\rho)$, 
$\rho\in{\mathbb R}\smallsetminus\{0\}$ & 
$\rho{\mathrm{prox}}_{\varphi/\rho^2}(x/\rho)$\\
\hline
\thenumperceuse\refstepcounter{numperceuse}\label{plct:propLeblon}  &
reflection & $\varphi(-x)$  & $-{\mathrm{prox}}_\varphi(-x)$\\
\hline
\thenumperceuse\refstepcounter{numperceuse} &

quadratic  \rule{0pt}{2.5ex} & 
$\varphi(x)+\alpha\|x\|^2/2+{u}^\top{x}+\gamma$
& ${\mathrm{prox}}_{\varphi/(\alpha+1)}
\big((x-u)/(\alpha+1)\big)$\\
& perturbation & $u\in{\mathbb R}^N$, $\alpha\geq 0$, 
$\gamma\in{\mathbb R}$ & \\
\hline
\thenumperceuse\refstepcounter{numperceuse}\label{plct:propA380}  &
conjugation & $\varphi^*(x)$ & $x-{\mathrm{prox}}_{\varphi}{x}$\\[1mm]
\hline
\thenumperceuse \refstepcounter{numperceuse}\label{plct:propvi} &
squared distance \rule{0pt}{2.5ex} & 
$\displaystyle{\frac12}d_C^2(x)$ & 
$\displaystyle{\frac12}(x+P_Cx)$\\[2mm]
\hline
&&&\\[-3mm]
\thenumperceuse \refstepcounter{numperceuse}\label{plct:propvii} &
Moreau envelope \rule{0pt}{2.5ex} & 
$\widetilde{\varphi}(x)=\displaystyle{\inf_{y\in{\mathbb R}^N}}
\varphi(y)+\frac12\|x-y\|^2$
& $\displaystyle{\frac12}\big(x+{\mathrm{prox}}_{2\varphi}x\big)$\\
\hline
&&&\\[-3mm]
\thenumperceuse \refstepcounter{numperceuse}\label{plct:propvII} &
Moreau complement \rule{0pt}{2.5ex} & 
$\displaystyle{\frac12}\|\cdot\|^2-\widetilde{\varphi}(x)$
& $x-{\mathrm{prox}}_{\varphi/2}(x/2)$\\[2mm]
\hline
&&&\\[-3mm]
\thenumperceuse \refstepcounter{numperceuse}\label{plct:prop8}  &
decomposition & 
$\sum_{k=1}^N\phi_k({x}^\top{b_k})$ & 
$\sum_{k=1}^N{\mathrm{prox}}_{\phi_k}\!({x}^\top{b_k})b_k$\\
&  $\begin{array}{l} 
\text{in an orthonormal}\\
\text{basis}\; (b_k)_{1\le k \le N}
\end{array}$
 & $\begin{array}{l}
\phi_k\in\Gamma_0({\mathbb R})\\
\end{array}$ & \\[3mm]
\hline
\thenumperceuse \refstepcounter{numperceuse} &
semi-orthogonal \rule{0pt}{2.5ex} &  $\varphi(Lx)$ &  $x+\nu^{-1}L^\top
\big({\mathrm{prox}}_{\nu \varphi}(Lx)-Lx\big)$\\
& linear transform &$L\in{\mathbb R}^{M\times N}$, 
$L L^\top = \nu I$, $\nu>0$ & \\
\hline
\thenumperceuse \refstepcounter{numperceuse} \label{plct:propviii} &
quadratic function \rule{0pt}{2.5ex} & $\gamma\|Lx-y\|^2/2$ & 
$(I+\gamma L^\top L)^{-1}(x+\gamma L^\top y)$\\
& & $L\in{\mathbb R}^{M\times N}$, $\gamma > 0$, $y\in{\mathbb R}^M$ & \\
\hline
\thenumperceuse \refstepcounter{numperceuse}  &
indicator function &
$\iota_C(x) = 
\begin{cases}
0 & \text{if}\;\;x\in C\\
+\infty & \text{otherwise}
\end{cases}$ 
& $P_C x$\\
\hline
\thenumperceuse\refstepcounter{numperceuse}  &
distance function & $\gamma d_C(x)$, $\gamma > 0$ &
$\begin{cases}
x+\gamma (P_Cx -x)/d_C(x)\\
\hspace*{10mm}\text{if}\;\;d_C(x)>\gamma\\
P_C x \hspace*{0.55cm} \text{otherwise}
\end{cases}
$\\
\hline
$\begin{array}{l}
\text{\thenumperceuse\refstepcounter{numperceuse} }\\
~
\end{array}$
&
$\begin{array}{l}
\text{function of}\\
\text{distance}
\end{array}$
& $\begin{array}{l}
\phi(d_C(x))\\
\phi\in\Gamma_0({\mathbb R})\;\text{even, differentiable}\\
\text{at $0$ with $\phi'(0)=0$}
\end{array}$
& 
$
\begin{cases}
x+ \left(1-
{\displaystyle{\frac{{\mathrm{prox}}_\phi d_C(x)}{d_C(x)}}}\right)
(P_C x-x)\\
\hspace*{22mm}\text{if}\;\;x\notin C\\
x \hspace*{20.8mm} \text{otherwise}
\end{cases}
$\\
\hline
 \thenumperceuse\refstepcounter{numperceuse} &
support function & $\sigma_C(x)$ & $x -P_Cx$\\
\hline
$\begin{array}{l}
\text{\thenumperceuse\refstepcounter{numperceuse}}\\
~ 
\end{array}$
&
$\begin{array}{l}
\text{thresholding}\\
~
\end{array}$
& $\begin{array}{l}
\sigma_C(x)+\phi(\|x\|)\\
\text{$\phi \in \Gamma_0({\mathbb R})$ even}\\
\text{and not constant}
\end{array}$
& $\begin{cases}
{\displaystyle{\frac{{\mathrm{prox}}_\phi d_C(x)}{d_C(x)}}}(x-P_C x)\\
\hspace*{11mm}\text{if}\;\;d_C(x) > 
\max\operatorname{Argmin} \phi\\
x-P_C x\hskip 2mm \text{otherwise}
\end{cases}$\\
\hline
\end{tabular}
\end{table}

\section{Forward-backward splitting}
\label{plcsec:3}

In this section, we consider the case of $m=2$ functions in 
\eqref{plce:1}, one of which is smooth.

\begin{problem}
\label{plcprob:2}
Let $f_1\in\Gamma_0({\mathbb R}^N)$, let 
$f_2\colon{\mathbb R}^N\to{\mathbb R}$ be convex 
and differentiable with a $\beta$-Lipschitz continuous 
gradient $\nabla f_2$, i.e., 
\begin{equation}
\label{plce:4}
(\forall (x,y)\in{\mathbb R}^N\times{\mathbb R}^N)\quad
\|\nabla f_2(x)-\nabla f_2(y)\|\leq\beta\|x-y\|,
\end{equation}
where $\beta\in\,\left]0,+\infty\right[$. 
Suppose that $f_1(x)+f_2(x)\to+\infty$ as 
$\|x\|\to+\infty$. The problem is to 
\begin{equation}
\label{plce:2}
\underset{x\in{\mathbb R}^N}{\text{minimize}}\;\;f_1(x)+f_2(x).
\end{equation}
\end{problem}

It can be shown \cite{plcSmms05} that Problem~\ref{plcprob:2} 
admits at least one solution and that, for any 
$\gamma\in\,\left]0,+\infty\right[$, its solutions 
are characterized by the fixed point equation
\begin{equation}
\label{plce:5}
x={\mathrm{prox}}_{\gamma f_1}\big(x-\gamma\nabla f_2(x)\big).
\end{equation}
This equation suggests the possibility of iterating 
\begin{equation}
\label{plce:24}
x_{n+1}=
\underbrace{{\mathrm{prox}}_{\gamma_n f_1}}_{\text{backward step}}
\big(\underbrace{x_n-\gamma_n\nabla f_2(x_n)}_{\text{forward step}}\big)
\end{equation}
for values of the step-size parameter $\gamma_n$ in a suitable 
bounded interval. This type of scheme is known as a
\emph{forward-backward} splitting algorithm for, using the
terminology used in discretization schemes in numerical analysis
\cite{plcVarg00}, it can be broken up into a forward (explicit) 
gradient step using the function $f_2$, and a backward (implicit) 
step using the function $f_1$. The forward-backward algorithm 
finds its roots in the projected gradient method \cite{plcLevi66} and 
in decomposition methods for solving variational inequalities 
\cite{plcMerc79,plcSibo70}. More recent forms of the algorithm and 
refinements can be found in 
\cite{plcBred09,plcChen97,plcOpti04,plcHale08,plcTsen91}.
Let us note that, on the one hand, when $f_1=0$, 
\eqref{plce:24} reduces to the \emph{gradient method}
\begin{equation}
\label{plce:25}
x_{n+1}=x_n-\gamma_n\nabla f_2(x_n)
\end{equation}
for minimizing a function with a Lipschitz continuous gradient
\cite{plcBert97,plcDunn76}. 
On the other hand, when $f_2=0$, \eqref{plce:24} reduces to the 
\emph{proximal point algorithm}
\begin{equation}
\label{plce:26}
x_{n+1}={\mathrm{prox}}_{\gamma_n f_1}x_n
\end{equation}
for minimizing a nondifferentiable function 
\cite{plcBrez78,plcOpti04,plcLema89,plcMart70,plcRock76}.
The forward-backward algorithm can therefore be considered as a
combination of these two basic schemes. The following version 
incorporates relaxation parameters $(\lambda_n)_{n\in{\mathbb N}}$.

\begin{algorithm}[Forward-backward algorithm]
\label{plcalgo:FB}
~\\
Fix $\varepsilon\in\:]0,\min\{1,1/\beta\}[$, $x_0\in{\mathbb R}^N$\\
For\;$n=0,1,\ldots$
\begin{equation}
\label{plce:FB}
\hskip -69mm
\begin{array}{l}
\left\lfloor
\begin{array}{l}
\gamma_n \in [\varepsilon,2/\beta-\varepsilon]\\
y_n=x_n-\gamma_n\nabla f_2(x_n)\\
\lambda_n\in [\varepsilon,1]\\
x_{n+1}=x_n+\lambda_n({\mathrm{prox}}_{\gamma_n f_1}y_n-x_n).
\end{array}
\right.\\
\end{array}
\end{equation}
\end{algorithm}

\begin{proposition}{\rm \cite{plcSmms05}}
\label{plcp:2}
Every sequence $(x_n)_{n\in{\mathbb N}}$ generated by 
Algorithm~\ref{plcalgo:FB} converges to a solution to 
Problem~\ref{plcprob:2}.
\end{proposition}

The above forward-backward algorithm features varying 
step-sizes $(\gamma_n)_{n\in{\mathbb N}}$ but its relaxation parameters
$(\lambda_n)_{n\in{\mathbb N}}$ cannot exceed 1. The following variant
uses constant step-sizes and larger relaxation parameters.

\begin{algorithm}[Constant-step forward-backward algorithm]
\label{plcalgo:FB2}
~\\
Fix $\varepsilon\in\:]0,3/4[$ and $x_0\in{\mathbb R}^N$\\
For\;$n=0,1,\ldots$
\begin{equation}
\hskip -6.7cm
\begin{array}{l}
\left\lfloor
\begin{array}{l}
y_n=x_n-\beta^{-1}\nabla f_2(x_n)\\
\lambda_n\in [\varepsilon,3/2-\varepsilon]\\
x_{n+1}=x_n+\lambda_n({\mathrm{prox}}_{\beta^{-1} f_1}y_n-x_n).
\end{array}
\right.\\
\end{array}
\end{equation}
\end{algorithm}

\begin{proposition}{\rm \cite{plcLivr10}}
\label{plcp:2''}
Every sequence $(x_n)_{n\in{\mathbb N}}$ generated by 
Algorithm~\ref{plcalgo:FB2} converges to a solution to 
Problem~\ref{plcprob:2}.
\end{proposition}

Although they may have limited impact on actual numerical performance, 
it may be of interest to know whether linear convergence rates are
available for the forward-backward algorithm.
In general, the answer is negative: even in the simple setting of 
Example~\ref{plcex:2} below, linear convergence of the iterates
$(x_n)_{n\in{\mathbb N}}$ generated by 
Algorithm~\ref{plcalgo:FB} fails \cite{plcBaus94,plcYoul87}.
Nonetheless it can be achieved at the expense of additional 
assumptions on the problem 
\cite{plcBaus96,plcBred08,plcChen97,plcDunn76,plcLema98,plcMerc79,%
plcMerc80,plcRock76,plcSibo70,plcZhuc95}.

Another type of convergence rate is that pertaining to the
objective values $(f_1(x_n)+f_2(x_n))_{n\in{\mathbb N}}$. 
This rate has been investigated in several places 
\cite{plcBec09a,plcBred08,plcGule91} and variants
of Algorithm~\ref{plcalgo:FB} have been developed to improve it
\cite{plcBec09a,plcBec09b,plcGule92,plcNest05,plcNest07,
plcTsen08,plcWeis09}
in the spirit of classical work by Nesterov \cite{plcNest83}.
It is important to note that the convergence of the 
sequence of iterates $(x_n)_{n\in{\mathbb N}}$, which is often crucial
in practice, is no longer guaranteed in general in such variants.
The proximal gradient method proposed in 
\cite{plcBec09a,plcBec09b} assumes the following form.

\begin{algorithm}[Beck-Teboulle proximal gradient algorithm]
\label{plcalgo:FISTA}
~\\
Fix $x_0\in{\mathbb R}^N$, set $z_0=x_0$ and $t_0=1$\\
For\;$n=0,1,\ldots$
\begin{equation}
\hskip -7.7cm
\begin{array}{l}
\left\lfloor
\begin{array}{l}
y_{n}=z_n-\beta^{-1}\nabla f_2(z_n)\\[2mm]
x_{n+1}={\mathrm{prox}}_{\beta^{-1}f_1}y_n\\[2mm]
\displaystyle{t_{n+1}=\frac{1+\sqrt{4t_n^2+1}}{2}}\\[2mm]
\displaystyle{\lambda_{n} = 1+\frac{t_n-1}{t_{n+1}}}\\[3mm]
\displaystyle{z_{n+1}=x_{n}+\lambda_n(x_{n+1}-x_n)}.
\end{array}
\right.\\
\end{array}
\end{equation}
\end{algorithm}

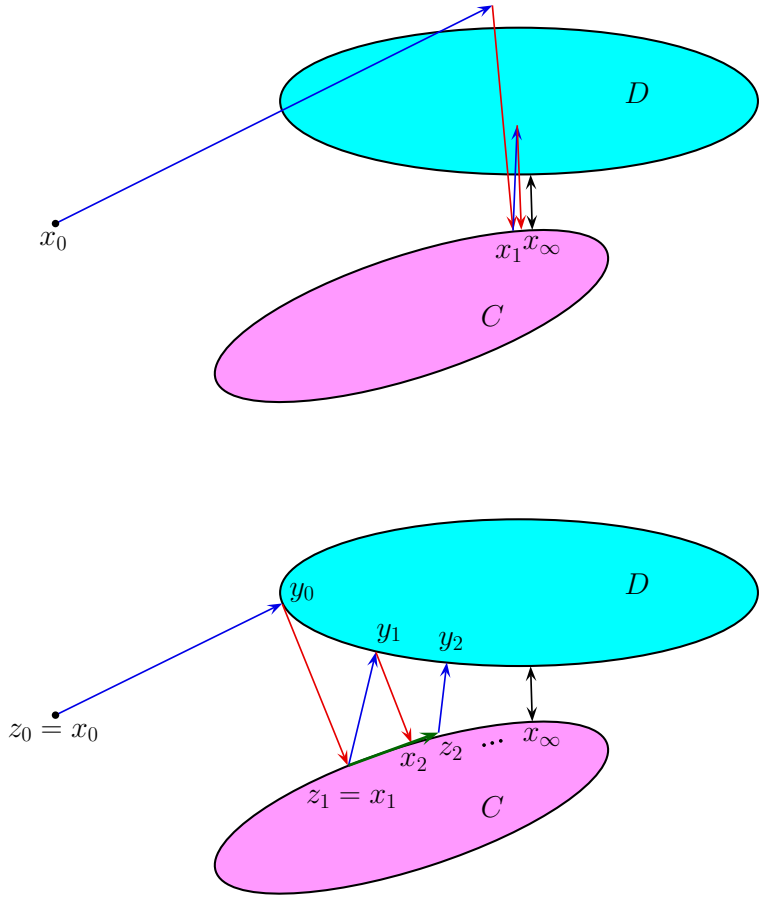
\begin{figure}
\begin{center}

\scalebox{0.55}{
\begin{pspicture}(0,-4.85)(18.24,6.98)
\definecolor{color96b}{rgb}{0.0,1.0,1.0}
\definecolor{color98b}{rgb}{0.0,0.0,0.9}
\definecolor{color95b}{rgb}{1.0,0.6,1.0}
\definecolor{color97b}{rgb}{0.9,0.0,0.0}
\rput{18.0}(-0.18953674,-2.8316424)%
{\psellipse[linewidth=0.05,dimen=outer,fillstyle=solid,fillcolor=color95b]%
(8.844375,-2.0141652)(5.0,1.5)}
\psellipse[linewidth=0.05,dimen=outer,fillstyle=solid,fillcolor=color96b]%
(11.444375,3.185835)(5.8,1.8)
\psline[linewidth=0.04cm,linecolor=color98b,arrowsize=0.18cm 2.0,%
arrowlength=1.4,arrowinset=0.4]{->}(0.204375,0.20583488)(10.8,5.5)
\psline[linewidth=0.04cm,linecolor=color97b,arrowsize=0.18cm 2.0,%
arrowlength=1.4,arrowinset=0.4]{->}(10.8,5.5)(11.3,0.05)
\psline[linewidth=0.04cm,linecolor=color98b,arrowsize=0.18cm 2.0,%
arrowlength=1.4,arrowinset=0.4]{->}(11.3,0.05)(11.4,2.6)
\psline[linewidth=0.04cm,linecolor=color97b,arrowsize=0.18cm 2.0,%
arrowlength=1.4,arrowinset=0.4]{->}(11.4,2.6)(11.5,0.07)
\psline[linewidth=0.04cm,arrowsize=0.18cm 2.0,arrowlength=1.4,%
arrowinset=0.4]{<->}(11.764375,0.06583487)(11.724375,1.4058349)

\psdots[dotsize=0.18](0.24,0.22)
\rput(10.8,-2.0){\huge $C$}
\rput(14.3,3.4){\huge $D$}
\rput(0.19,-0.20){\huge $x_0$}
\rput(11.2,-0.5){\huge $x_1$}
\rput(12.0,-0.3){\huge $x_\infty$}
\end{pspicture} 
}

\scalebox{0.55}{
\begin{pspicture}(0,-4.85)(18.24,6.98)
\definecolor{color96b}{rgb}{0.0,1.0,1.0}
\definecolor{dgreen}{rgb}{0.0,0.4,0.0}
\definecolor{color98b}{rgb}{0.0,0.0,0.9}
\definecolor{color95b}{rgb}{1.0,0.6,1.0}
\definecolor{color97b}{rgb}{0.9,0.0,0.0}
\rput{18.0}(-0.18953674,-2.8316424)%
{\psellipse[linewidth=0.05,dimen=outer,fillstyle=solid,fillcolor=color95b]%
(8.844375,-2.0141652)(5.0,1.5)}
\psellipse[linewidth=0.05,dimen=outer,fillstyle=solid,fillcolor=color96b]%
(11.444375,3.185835)(5.8,1.8)
\psline[linewidth=0.04cm,linecolor=color98b,arrowsize=0.18cm 2.0,%
arrowlength=1.4,arrowinset=0.4]{->}(0.204375,0.20583488)(5.724375,2.925835)
\psline[linewidth=0.04cm,linecolor=color97b,arrowsize=0.18cm 2.0,%
arrowlength=1.4,arrowinset=0.4]{->}(5.7243,2.9058)(7.304375,-0.994165)
\psline[linewidth=0.04cm,linecolor=color98b,%
arrowsize=0.18cm 2.0,arrowlength=1.4,%
arrowinset=0.4]{->}(7.324375,-0.9941651)(7.984375,1.7458348)
\psline[linewidth=0.04cm,linecolor=color97b,%
arrowsize=0.18cm 2.0,arrowlength=1.4,%
arrowinset=0.4]{->}(7.984375,1.7458348)(8.844375,-0.45416513)
\psline[linewidth=0.07cm,linecolor=dgreen,%
arrowsize=0.18cm 2.0,arrowlength=1.4,%
arrowinset=0.4]{->}(7.324375,-0.9941651)(9.5,-0.2)
\psline[linewidth=0.04cm,linecolor=color98b,%
arrowsize=0.18cm 2.0,arrowlength=1.4,%
arrowinset=0.4]{->}(9.5,-0.2)(9.7,1.5)

\psline[linewidth=0.04cm,arrowsize=0.18cm 2.0,arrowlength=1.4,%
arrowinset=0.4]{<->}(11.764375,0.06583487)(11.724375,1.4058349)
\psdots[dotsize=0.18](0.24,0.22)
\rput(10.8,-2.0){\huge $C$}
\rput(14.3,3.4){\huge $D$}
\rput(0.19,-0.20){\huge $z_0=x_0$}
\rput(7.4,-1.8){\huge $z_1=x_1$}
\rput(6.2,3.2){\huge $y_0$}
\rput(8.3,2.2){\huge $y_1$}
\rput(9.8,2.0){\huge $y_2$}
\rput(8.9,-0.91){\huge $x_2$}
\rput(9.8,-0.65){\huge $z_2$}
\rput(12.0,-0.3){\huge $x_\infty$}
\psdots[dotsize=0.09](10.6,-0.50)
\psdots[dotsize=0.09](10.8,-0.45)
\psdots[dotsize=0.09](11.0,-0.40)
\end{pspicture} 
}

\caption{Forward-backward versus Beck-Teboulle : 
As in Example~\ref{plcex:2}, let $C$ and $D$ be two closed convex sets 
and consider the problem \eqref{plce:11} of finding a point $x_\infty$
in $C$ at minimum distance from $D$. Let us set $f_1=\iota_C$ and 
$f_2=d_D^2/2$. 
Top: The forward-backward algorithm 
with $\gamma_n\equiv 1.9$ and $\lambda_n\equiv 1$.
As seen in Example~\ref{plcex:2}, it reduces to the 
alternating projection method \eqref{plce:14}. 
Bottom: The Beck-Teboulle algorithm.}
\label{plcfig:e8}
\end{center}
\end{figure}

While little is known about the actual convergence of sequences
produced by Algorithm~\ref{plcalgo:FISTA},
the $O(1/n^2)$ rate of convergence of the objective function 
they achieve is optimal \cite{plcNemi83}, although the practical
impact of such property is not always manifest in concrete 
problems (see Figure~\ref{plcfig:e8} for a comparison with the 
Forward-Backward algorithm).

\begin{proposition}{\rm \cite{plcBec09a}}
\label{plcp:2'}
Assume that, for every $y\in\text{dom}\,f_1$, 
$\partial f_1(y)\neq\varnothing$, and
let $x$ be a solution to Problem~\ref{plcprob:2}.
Then every sequence $(x_n)_{n\in{\mathbb N}}$ generated by 
Algorithm~\ref{plcalgo:FISTA} satisfies
\begin{equation}
(\forall n\in {\mathbb N}\smallsetminus\{0\})\quad
f_1(x_n)+f_2(x_n)\leq f_1({x})+f_2({x}) 
+\frac{2\beta\|x_0-{x}\|^2}{(n+1)^2}.
\end{equation}
\end{proposition}

Other variations of the forward-backward algorithm
have also been reported to yield improved convergence profiles 
\cite{plcBiou07,plcPeyr09,plcLori09,plcVone08,plcVone09}.

Problem~\ref{plcprob:2} and Proposition~\ref{plcp:2} 
cover a wide variety of signal processing problems and 
solution methods \cite{plcSmms05}. 
For the sake of illustration, let us provide a few examples.
For notational convenience, we set $\lambda_n\equiv 1$ in 
Algorithm~\ref{plcalgo:FB}, which reduces the updating 
rule to \eqref{plce:24}.

\refstepcounter{numscie}
\begin{table}[hptb]
\centering
\caption{Proximity operator of 
$\phi\in\Gamma_0({\mathbb R})$; $\alpha\in{\mathbb R}$,
$\kappa>0$, $\underline{\kappa}>0$, $\overline{\kappa}>0$, 
$\omega>0$, $\underline{\omega}<\overline{\omega}$, $q>1$, 
$\tau\geq 0$ 
\cite{plcChau07,plcSiop07,plcSmms05}.}
\label{plct:real}
{\small
\begin{tabular}{|p{0.4cm} l|l|}
\hline
& $\phi(x)$ & ${\mathrm{prox}}_{\phi} x$\\
\hline 
\hline
\thenumscie\refstepcounter{numscie}\label{plct:propsoft} &
$\iota_{[\underline{\omega},\overline{\omega}]}(x)$ & 
$P_{[\underline{\omega},\overline{\omega}]}\,x$\\[1mm]
\hline 
\thenumscie\refstepcounter{numscie} &
$\sigma_{[\underline{\omega},\overline{\omega}]}(x)=\begin{cases}
\underline{\omega}x&\text{if}\;\;x<0\\
0&\text{if}\;\;x=0\\
\overline{\omega}x&\text{otherwise}
\end{cases}$ & 
${\operatorname{soft}}_{[\underline{\omega},\overline{\omega}]}
(x) =\begin{cases}
x-\underline{\omega}&\text{if}\;\;x<\underline{\omega}\\
0&\text{if}\;\;x\in[\underline{\omega},\overline{\omega}]\\
x-\overline{\omega}&\text{if}\;\;x>\overline{\omega}
\end{cases}$\\
\hline
\thenumscie\refstepcounter{numscie} &
$\begin{array}{l}
\psi(x)+\sigma_{[\underline{\omega},\overline{\omega}]}(x)\\
\psi\in \Gamma_0({\mathbb R})\; \text{differentiable at $0$}\\
\psi'(0) = 0
\end{array}$
& ${\mathrm{prox}}_\psi\big(
{\operatorname{soft}}_{[\underline{\omega},\overline{\omega}]}
(x)\big)$\\
\hline
\thenumscie\refstepcounter{numscie} &
$\max\{|x|-\omega,0\}$ & 
$\begin{cases}
x & \text{if $|x|<\omega$ }\\
\operatorname{sign}(x)\omega &\text{if $\omega\le |x|\le 2\omega$}\\ 
\operatorname{sign}(x) (|x| -\omega) & \text{if $|x| > 2\omega$}
\end{cases}$\\
\hline
\thenumscie\refstepcounter{numscie} &
$\kappa|x|^{q}$ & $\begin{array}{l}
\operatorname{sign}(x)p,\\
\text{where}\;
p \ge 0\;\text{and}\;p+q\kappa p^{q-1}=|x|
\end{array}$
\\
\hline
\thenumscie\refstepcounter{numscie} &
$\begin{cases}
\kappa x^2 & \text{if}\;\;|x|\leq\omega/\sqrt{2\kappa}\\[1mm]
\omega\sqrt{2\kappa}|x|-
{\omega^2}/{2} & \text{otherwise}
\end{cases}$ &
 $\begin{cases}
x/(2\kappa+1) 
&\text{if}\;\;|x|\leq\omega(2\kappa+1)/\sqrt{2\kappa}\\[1mm]
x-\omega\sqrt{2\kappa}\operatorname{sign}(x)
&\text{otherwise}
\end{cases}$\\
\hline
\thenumscie\refstepcounter{numscie} &
$\omega|x|+\tau|x|^2+\kappa|x|^q$ & 
$\displaystyle\operatorname{sign}(x)
{\mathrm{prox}}_{\kappa|\cdot|^q/(2\tau+1)}
\frac{\max\{|x|-\omega,0\}}{2\tau+1}$\\[2mm]
\hline
\thenumscie\refstepcounter{numscie} &
$\omega|x|-\ln(1+\omega|x|)$ & 
$\begin{array}{l}
(2\omega)^{-1}\operatorname{sign}(x)
\Big(\omega|x|-\omega^2-1\\
\hspace*{1.9cm}+
\sqrt{\big|\omega|x|-\omega^2-1 \big|^2+4\omega|x|}\Big)
\end{array}$\\
\hline
\thenumscie\refstepcounter{numscie} &
$\begin{cases}
\omega x&\text{if}\;\;x\geq 0\\
+\infty & \text{otherwise}
\end{cases}$ & $\begin{cases}
x-\omega&\text{if}\;\;x\geq\omega\\
0&\text{otherwise}
\end{cases}$\\
\hline
\thenumscie\refstepcounter{numscie} &$\begin{cases}
-\omega x^{1/q}&\text{if}\;\;x\ge 0\\
+\infty & \text{otherwise}
\end{cases}$ & $\begin{array}{l}
p^{1/q},\\
\text{where}\;p > 0\;\text{and}\;
p^{2q-1} -x p^{q-1} = q^{-1}\omega
\end{array}$\\
\hline
\thenumscie\refstepcounter{numscie} &
$\begin{cases}
\omega x^{-q}&\text{if}\;\;x > 0\\
+\infty & \text{otherwise}
\end{cases}$ & $\begin{array}{l}
p > 0\\
\text{such that}\; p^{q+2} -x p^{q+1}=\omega q
\end{array}$\\
\hline 
\thenumscie\refstepcounter{numscie}\label{plct:propfig} &
$\begin{cases}
x\ln(x)&\text{if}\;\;x>0\\
0&\text{if}\;\;x=0\\
+\infty&\text{otherwise}
\end{cases}$ & 
$\begin{array}{l}
W(e^{x-1}),\\
\text{where $W$ is the Lambert W-function}\\
\end{array}$\\
\hline
\thenumscie\refstepcounter{numscie} &
$\begin{cases}
-\ln(x -\underline{\omega} )+\ln(-\underline{\omega})
& \text{if}\;\; x\in\left]\underline{\omega},0\right]\\
-\ln(\overline{\omega}-x)+\ln(\overline{\omega})&
\text{if}\;\;x\in\left]0,\overline{\omega}\right[\\
+\infty &\text{otherwise}
\end{cases}$ \rule{0pt}{7ex} &   $\begin{cases}
\displaystyle{\frac12\Big(x+\underline{\omega}+
\sqrt{|x-\underline{\omega}|^2+4}\Big)}
& \text{if}\;\;x<1/{\underline\omega}\\[2mm]
\displaystyle{\frac12\Big(x+\overline{\omega}-
\sqrt{|x-\overline{\omega}|^2+4}\Big)} 
& \text{if}\;\;x>1/\overline{\omega}\\
0 & \text{otherwise}
\end{cases}$\\
& $\underline{\omega} < 0< \overline{\omega}$ &
\text{(see Figure~\ref{fig:6fevrier2010})}\\
\hline
\thenumscie\refstepcounter{numscie} &
$\begin{cases}
\displaystyle -\kappa\ln(x)+{\tau x^2}/{2}+\alpha x& \text{if}\;\; x>0\\
+\infty & \text{otherwise}
\end{cases}$ & $\displaystyle{\frac{1}{2(1+\tau)}}
\Big(x-\alpha+\sqrt{|x-\alpha|^2+4\kappa(1+\tau)}\Big)$\\
\hline
\thenumscie\refstepcounter{numscie} &
$\begin{cases}
-\kappa\ln(x)+\alpha x+\omega x^{-1} &\text{if}\;\;x>0\\
+\infty & \text{otherwise}
\end{cases}$ &$\begin{array}{l}
p > 0\\ 
\text{such that}\;p^3+(\alpha-x)p^2-\kappa p=\omega
\end{array}$\\
\hline
\thenumscie\refstepcounter{numscie}  &
$\begin{cases}
-\kappa\ln(x)+\omega x^q&\text{if}\;\;x>0\\
+\infty & \text{otherwise}
\end{cases}$ & $\begin{array}{l}
p>0\\ 
\text{such that}\;q\omega p^q+p^2-xp=\kappa
\end{array}$\\
\hline
\thenumscie\refstepcounter{numscie} &
$\begin{cases}
-\underline{\kappa} \ln(x-\underline{\omega})-\overline{\kappa}
\ln(\overline{\omega}-x)\\ 
\hspace*{2.9cm}\text{if}\;\;x 
\in\,]\underline{\omega},\overline{\omega}[\\
+\infty \hspace*{2.45cm}\text{otherwise}
\end{cases}$
& $\begin{array}{l}
p \in\, ]\underline{\omega},\overline{\omega}[\; \\
\text{such that}\;p^3-(\underline{\omega}+\overline{\omega}+x)p^2+
\\ \big(\underline{\omega}\overline{\omega}-\underline{\kappa}
-\overline{\kappa}+(\underline{\omega}+\overline{\omega})x\big)p=
\underline{\omega}\overline{\omega}x-\underline{\omega}
\overline{\kappa}
-\overline{\omega}\underline{\kappa}
\end{array}$
\\
\hline
\end{tabular}
}
\end{table}

\begin{figure}
\begin{center}
\scalebox{0.65} 
{
\begin{pspicture}(-9,-6.5)(10,8) 
\psline[linewidth=0.05cm,arrowsize=0.4cm,linestyle=solid]{->}(-8,0)(8,0)
\psline[linewidth=0.05cm,arrowsize=0.4cm,linestyle=solid]{->}(0,-6)(0,6)
\psplot[plotpoints=800,algebraic,linewidth=0.08cm,linecolor=blue]%
{-8.0}{-1.0}{5*(x-1+sqrt((x+1)^2+4))/2}
\psplot[plotpoints=800,algebraic,linewidth=0.08cm,linecolor=blue]%
{-1.0}{1.0}{0}
\psplot[plotpoints=800,algebraic,linewidth=0.08cm,linecolor=blue]%
{1.0}{8.0}{5*(x+1-sqrt((x-1)^2+4))/2}
\psline[linewidth=0.05cm,linestyle=dashed]{-}(0.0,5.0)(8.0,5.0)
\psline[linewidth=0.05cm,linestyle=dashed]{-}(-8.0,-5.0)(0.0,-5.0)
\rput(8.4,0){\Large${\xi}$}
\rput(0,6.4){\Large${{\mathrm{prox}}_\phi\,\xi}$}
\rput(1.0,-0.4){\Large${1/\overline{\omega}}$}
\rput(-1.0,0.4){\Large${1/\underline{\omega}}$}
\rput(-0.4,5.0){\Large${\overline{\omega}}$}
\rput(0.4,-5.0){\Large${\underline{\omega}}$}
\end{pspicture} 
}
\caption{Proximity operator of the function}
\[
\phi\colon{\mathbb R}\to\,\left]-\infty,+\infty\right]\colon\xi\mapsto
\begin{cases}
-\ln(\xi-\underline{\omega} )+\ln(-\underline{\omega})
& \text{if}\;\; \xi\in\left]\underline{\omega},0\right]\\
-\ln(\overline{\omega}-\xi)+\ln(\overline{\omega})&
\text{if}\;\;\xi\in\left]0,\overline{\omega}\right[\\
+\infty &\text{otherwise.}
\end{cases}
\]
\begin{minipage}{\textwidth}
The proximity operator thresholds over the interval
$[1/\underline{\omega},1/\overline{\omega}]$, and saturates
at $-\infty$ and $+\infty$ with asymptotes at $\underline{\omega}$
and $\overline{\omega}$, respectively 
(see Table~\ref{plct:real}.\ref{plct:propfig}
and \cite{plcSiop07}).
\end{minipage}
\label{fig:6fevrier2010}
\end{center}
\end{figure}
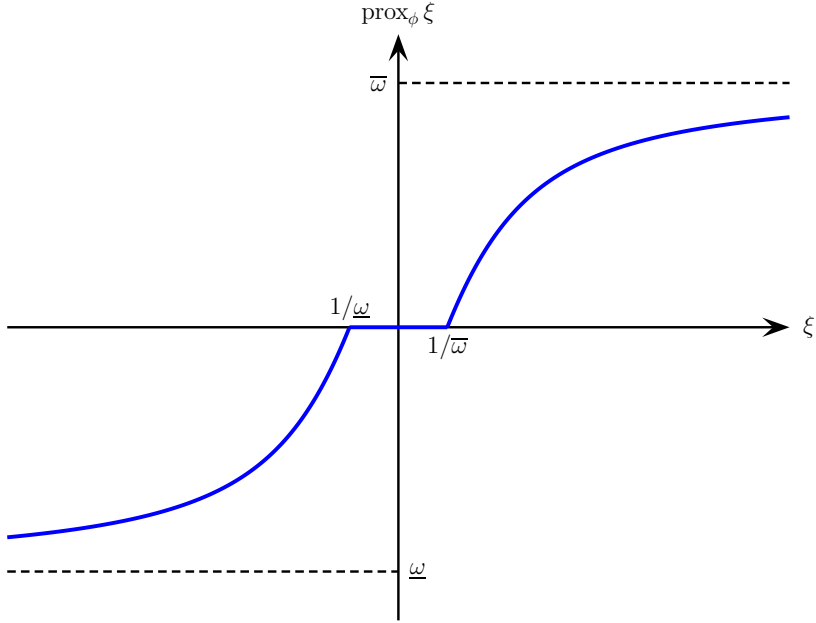

\begin{plcexample}[projected gradient]
\label{plcex:111}
In Problem~\ref{plcprob:2}, suppose that $f_1=\iota_{C}$,
where $C$ is a closed convex subset of ${\mathbb R}^N$ such that
$\{x\in C\:|\:f_2(x)\leq\eta\}$ is nonempty and bounded
for some $\eta\in{\mathbb R}$. Then we obtain
the constrained minimization problem
\begin{equation}
\label{plce:112}
\underset{x\in C}{\text{minimize}}\;\;f_2(x).
\end{equation}
Since ${\mathrm{prox}}_{\gamma f_1}=P_C$ 
(see Table~\ref{plct:prop}.\ref{plct:propviii}), 
the forward-backward iteration reduces to the 
\emph{projected gradient method}
\begin{equation}
\label{plce:113}
x_{n+1}=P_C\big(x_n-\gamma_n\nabla f_2(x_n)\big),
\;\;\varepsilon\le \gamma_n \le 2/\beta-\varepsilon.
\end{equation}
This algorithm has been used in numerous signal processing problems,
in particular in total variation denoising \cite{plcCham05}, in
image deblurring \cite{plcBenv10}, in pulse shape design 
\cite{plcSign99}, and in compressed sensing \cite{plcFigu07}.
\end{plcexample}

\begin{plcexample}[projected Landweber]
\label{plcex:1}
In Example~\ref{plcex:111}, setting $f_2\colon x\mapsto\|Lx-y\|^2/2$,
where $L\in{\mathbb R}^{M\times N}\smallsetminus\{0\}$ and 
$y\in{\mathbb R}^M$, yields 
the constrained least-squares problem
\begin{equation}
\label{plce:12}
\underset{x\in C}{\text{minimize}}\;\;\frac12\|Lx-y\|^2.
\end{equation}
Since $\nabla f_2\colon x\mapsto L^\top(Lx-y)$ has Lipschitz constant
$\beta=\|L\|^2$, \eqref{plce:113} yields the 
\emph{projected Landweber method} \cite{plcEick92} 
\begin{equation}
\label{plce:13}
x_{n+1}=P_C\big(x_n+\gamma_nL^\top(y-Lx_n)\big),\;\;
\varepsilon\leq\gamma_n\leq 2/\|L\|^2-\varepsilon.
\end{equation}
This method has been used in particular in computer vision 
\cite{plcJoha06} and in signal restoration \cite{plcTrus85}.
\end{plcexample}

\begin{plcexample}[backward-backward algorithm]
\label{plcex:6}
Let $f$ and $g$ be functions in $\Gamma_0({\mathbb R}^N)$.
Consider the problem 
\begin{equation}
\label{plce:17}
\underset{x\in{\mathbb R}^N}{\text{minimize}}\;\;f(x)+\widetilde{g}(x),
\end{equation}
where $\widetilde{g}$ is the Moreau envelope of $g$ (see
Table~\ref{plct:prop}.\ref{plct:propvi}), and suppose that 
$f(x)+\widetilde{g}(x)\to+\infty$ as $\|x\|\to+\infty$.
This is a special case of Problem~\ref{plcprob:2} with
$f_1=f$ and $f_2=\widetilde{g}$. Since
$\nabla f_2\colon x\mapsto x-{\mathrm{prox}}_g x$ has Lipschitz constant 
$\beta=1$ \cite{plcSmms05,plcMore65}, Proposition~\ref{plcp:2} with 
$\gamma_n\equiv 1$ asserts that the sequence $(x_n)_{n\in{\mathbb N}}$ 
generated by the \emph{backward-backward algorithm}
\begin{equation}
\label{plce:18}
x_{n+1}={\mathrm{prox}}_f({\mathrm{prox}}_gx_n)
\end{equation}
converges to a solution to \eqref{plce:17}. Detailed analyses of this
scheme can be found in \cite{plcAcke80,plcReic05,plcOpti04,plcPass79}.
\end{plcexample}

\begin{plcexample}[alternating projections]
\label{plcex:2}
In Example~\ref{plcex:6}, let $f$ and $g$ be respectively the indicator
functions of nonempty closed convex sets $C$ and $D$, 
one of which is bounded. Then \eqref{plce:17} amounts to
finding a signal $x$ in $C$ at closest distance from $D$, i.e., 
\begin{equation}
\label{plce:11}
\underset{x\in C}{\text{minimize}}\;\;\frac12d^2_D(x).
\end{equation}
Moreover, since ${\mathrm{prox}}_f=P_C$ and ${\mathrm{prox}}_g=P_D$, 
\eqref{plce:18}
yields the \emph{alternating projection method} 
\begin{equation}
\label{plce:14}
x_{n+1}=P_C(P_Dx_n),
\end{equation}
which was first analyzed in this context in \cite{plcChen59}.
Signal processing applications can be found in the areas
of spectral estimation \cite{plcGold85}, 
pulse shape design \cite{plcNoba95},
wavelet construction \cite{plcPesq96}, and signal synthesis
\cite{plcYoul86}.
\end{plcexample}

\begin{plcexample}[iterative thresholding]
\label{plcex:3}
Let $(b_k)_{1\leq k\leq N}$ be an orthonormal basis of ${\mathbb R}^N$,
let $(\omega_k)_{1\leq k\leq N}$ be strictly positive real numbers,
let $L\in{\mathbb R}^{M\times N}\smallsetminus\{0\}$, and let 
$y\in{\mathbb R}^M$. Consider the $\ell^1$--$\ell^2$ problem
\begin{equation}
\label{plce:15}
\underset{x\in{\mathbb R}^N}{\text{minimize}}\;\;\sum_{k=1}^{N}\omega_k
|{x}^\top{b_k}|+\frac12\|Lx-y\|^2.
\end{equation}
This type of formulation arises in signal recovery problems in which
$y$ is the observed signal and the original signal is known to have 
a sparse representation in the basis $(b_k)_{1\leq k\leq N}$, e.g.,
\cite{plcBect04,plcBiou07,plcDaub04,plcDemo02,
plcFigu03,plcFigu07,plcTibs96,plcTrop06}. 
We observe that \eqref{plce:15} is a special case of \eqref{plce:2} with 
\begin{equation}
\begin{cases}
f_1\colon x\mapsto\sum_{1\leq k\leq N}\omega_k|{x}^\top{b_k}|\\
f_2\colon x\mapsto\|Lx-y\|^2/2. 
\end{cases}
\end{equation}
Since ${\mathrm{prox}}_{\gamma f_1}\colon x\mapsto\sum_{1\leq k\leq N}
{\operatorname{soft}}_{[-\gamma\omega_k,\gamma\omega_k]}
({x}^\top{b_k})\:b_k$ (see
Table~\ref{plct:prop}.\ref{plct:propvii} 
and Table~\ref{plct:real}.\ref{plct:propsoft}), 
it follows from Proposition~\ref{plcp:2} 
that the sequence $(x_n)_{n\in{\mathbb N}}$ generated by the 
\emph{iterative thresholding algorithm}
\begin{multline}
\label{plce:16}
x_{n+1}=\sum_{k=1}^{N}\xi_{k,n}b_k,\quad\text{where}\\
\begin{cases}
\xi_{k,n}={\operatorname{soft}}_{[-\gamma_n\omega_k,\gamma_n\omega_k]}
{\big(x_n+\gamma_nL^\top(y-Lx_n)\big)}^\top{b_k}\\
\varepsilon\leq\gamma_n\leq2/\|L\|^2-\varepsilon,
\end{cases}
\end{multline}
converges to a solution to \eqref{plce:15}.
\end{plcexample}

Additional applications of the forward-backward algorithm in 
signal and image processing can be found in 
\cite{plcCaij08,plcCaij09,plcChan09,plcChaa09,plcRaym08,%
plcChau07,plcSiop07,plcSmms05,plcDaub07,plcForn07}.

\section{Douglas-Rachford splitting}
\label{plcsec:4}

The forward-backward algorithm of Section~\ref{plcsec:3} requires that
one of the functions be differentiable, with a Lipschitz continuous 
gradient. In this section, we relax this assumption.

\begin{problem}
\label{plcprob:3}
Let $f_1$ and $f_2$ be functions in $\Gamma_0({\mathbb R}^N)$ such that
\begin{equation}
\label{plce:cq1}
({\mathrm{ri}}\,{\mathrm{dom}}\, f_1)\cap
({\mathrm{ri}}\,{\mathrm{dom}}\, f_2)\neq\varnothing
\end{equation}
and $f_1(x)+f_2(x)\to+\infty$ as $\|x\|\to+\infty$. The problem is to
\begin{equation}
\label{plce:3}
\underset{x\in{\mathbb R}^N}{\text{minimize}}\;\;f_1(x)+f_2(x).
\end{equation}
\end{problem}

What is nowadays referred to as the \emph{Douglas-Rachford 
algorithm} goes back to a method originally proposed in \cite{plcDoug56} 
for solving matrix equations of the form $u=Ax+Bx$, where $A$ and $B$ 
are positive-definite matrices (see also \cite{plcVarg00}). 
The method was transformed in \cite{plcLieu69} to handle
nonlinear problems and further improved in \cite{plcLion79} to
address monotone inclusion problems. For further developments,
see  \cite{plcOpti04,plcJoca09,plcEcks92}.

Problem~\ref{plcprob:3} admits at least one solution and, 
for any $\gamma\in\,\left]0,+\infty\right[$, its solutions are
characterized by the two-level condition \cite{plcJsts07}
\begin{equation}
\label{plce:40}
\begin{cases}
x={\mathrm{prox}}_{\gamma f_2}y\\
{\mathrm{prox}}_{\gamma f_2}y=
{\mathrm{prox}}_{\gamma f_1}(2{\mathrm{prox}}_{\gamma f_2}y-y),
\end{cases}
\end{equation}
which motivates the following scheme.

\begin{algorithm}[Douglas-Rachford algorithm]
\label{plcalgo:DR}
~\\
Fix $\varepsilon\in\:]0,1[$,  $\gamma>0$, $y_0\in{\mathbb R}^N$\\
For\;$n=0,1,\ldots$
\begin{equation}
\hskip -5.8cm
\begin{array}{l}
\left\lfloor
\begin{array}{l}
x_n={\mathrm{prox}}_{\gamma f_2}y_n\\
\lambda_n\in [\varepsilon,2-\varepsilon]\\
y_{n+1}=y_n+\lambda_n
\big({\mathrm{prox}}_{\gamma f_1}\big(2x_n-y_n\big)-x_n\big). 
\end{array}
\right.
\end{array}
\end{equation}
\end{algorithm}

\begin{proposition}{\rm\cite{plcJsts07}}
\label{plcp:3}
Every sequence $(x_n)_{n\in{\mathbb N}}$ generated by 
Algorithm~\ref{plcalgo:DR} converges to a solution 
to Problem~\ref{plcprob:3}.
\end{proposition}

Just like the forward-backward algorithm, the Douglas-Rachford algorithm 
operates by splitting since it employs the functions $f_1$ 
and $f_2$ separately. It can be viewed as more general in scope
than the forward-backward algorithm in that it does not require 
that any of the functions have a Lipschitz continuous gradient. 
However, this observation must be weighed against the fact that
it may be more demanding numerically as it requires the 
implementation of two proximal steps at each iteration,
whereas only one is needed in the forward-backward algorithm.
In some problems, both may be easily implementable
(see Fig.~\ref{plcfig:A380} for an example) and it is not
clear a priori which algorithm may be more efficient.

\begin{figure}
\begin{center}
\scalebox{0.55}{
\begin{pspicture}(0,-4.85)(18.24,6.98)
\definecolor{color96b}{rgb}{0.0,1.0,1.0}
\definecolor{color98b}{rgb}{0.0,0.0,0.9}
\definecolor{color95b}{rgb}{1.0,0.6,1.0}
\definecolor{color97b}{rgb}{0.9,0.0,0.0}
\rput{18.0}(-0.18953674,-2.8316424)%
{\psellipse[linewidth=0.05,dimen=outer,fillstyle=solid,fillcolor=color95b]%
(8.844375,-2.0141652)(5.0,1.5)}
\psellipse[linewidth=0.05,dimen=outer,fillstyle=solid,fillcolor=color96b]%
(11.444375,3.185835)(5.8,1.8)
\psline[linewidth=0.04cm,linecolor=color98b,arrowsize=0.18cm 2.0,%
arrowlength=1.4,arrowinset=0.4]{->}(0.204375,0.20583488)(5.724375,2.925835)
\psline[linewidth=0.04cm,linecolor=color97b,arrowsize=0.18cm 2.0,%
arrowlength=1.4,arrowinset=0.4]{->}(5.7243,2.9058)(7.304375,-0.994165)
\psline[linewidth=0.04cm,linecolor=color98b,%
arrowsize=0.18cm 2.0,arrowlength=1.4,%
arrowinset=0.4]{->}(7.324375,-0.9941651)(7.984375,1.7458348)
\psline[linewidth=0.04cm,linecolor=color97b,%
arrowsize=0.18cm 2.0,arrowlength=1.4,%
arrowinset=0.4]{->}(7.984375,1.7458348)(8.844375,-0.45416513)
\psline[linewidth=0.04cm,linecolor=color98b,%
arrowsize=0.18cm 2.0,arrowlength=1.4,%
arrowinset=0.4]{->}(8.844375,-0.45416513)(9.184375,1.5658349)
\psline[linewidth=0.04cm,linecolor=color97b,%
arrowsize=0.18cm 2.0,arrowlength=1.4,%
arrowinset=0.4]{->}(9.184375,1.5458349)(9.744375,-0.23416513)
\psline[linewidth=0.04cm,linecolor=color98b,%
arrowsize=0.18cm 2.0,arrowlength=1.4,arrowinset=0.4]{->}%
(9.724375,-0.2541651)(10.04,1.42)
\psline[linewidth=0.04cm,arrowsize=0.18cm 2.0,arrowlength=1.4,%
arrowinset=0.4]{<->}(11.764375,0.06583487)(11.724375,1.4058349)
\psdots[dotsize=0.18](0.24,0.22)
\rput(10.8,-2.0){\huge $C$}
\rput(14.3,3.4){\huge $D$}
\rput(0.19,-0.20){\huge $x_0$}
\rput(7.4,-1.5){\huge $x_1$}
\rput(8.9,-0.91){\huge $x_2$}
\rput(9.9,-0.65){\huge $x_3$}
\rput(12.0,-0.3){\huge $x_\infty$}
\psdots[dotsize=0.09](10.6,-0.50)
\psdots[dotsize=0.09](10.8,-0.45)
\psdots[dotsize=0.09](11.0,-0.40)
\end{pspicture} 
}

\vskip 10mm

\scalebox{0.55}{
\begin{pspicture}(0,-5.35)(17.24,4.98)
\definecolor{color96b}{rgb}{0.0,1.0,1.0}
\definecolor{color98b}{rgb}{0.0,0.0,0.9}
\definecolor{color95b}{rgb}{1.0,0.6,1.0}
\definecolor{color97b}{rgb}{0.9,0.0,0.0}
\rput{18.0}(-0.18953674,-2.8316424)%
{\psellipse[linewidth=0.05,dimen=outer,fillstyle=solid,fillcolor=color95b]%
(8.844375,-2.0141652)(5.0,1.5)}
\psellipse[linewidth=0.05,dimen=outer,fillstyle=solid,fillcolor=color96b]%
(11.444375,3.185835)(5.8,1.8)
\psline[linewidth=0.04cm,linecolor=color98b,arrowsize=0.18cm 2.0,%
arrowlength=1.4,arrowinset=0.4]{->}(0.204375,0.20583488)(5.724375,2.925835)
\psline[linewidth=0.04cm,linecolor=color97b,arrowsize=0.18cm 2.0,%
arrowlength=1.4,arrowinset=0.4]{->}(5.7243,2.9058)(7.304375,-0.994165)
\psline[linewidth=0.04cm,arrowsize=0.18cm 2.0,arrowlength=1.4,%
arrowinset=0.4]{<->}(11.764375,0.06583487)(11.724375,1.4058349)
\psline[linewidth=0.04cm,arrowsize=0.18cm 2.0,arrowlength=1.4,%
arrowinset=0.4]{->}(7.3371873,-0.97416514)(4.5771875,-2.334165)
\psline[linewidth=0.04cm,linecolor=color98b,arrowsize=0.18cm 2.0,%
arrowlength=1.4,arrowinset=0.4]{->}(4.5971,-2.334)(6.53718,2.2458)
\psline[linewidth=0.04cm,linecolor=color97b,arrowsize=0.18cm 2.0,%
arrowlength=1.4,arrowinset=0.4]{->}(6.537,2.2458)(7.6371,-0.85416)
\psline[linewidth=0.04cm,arrowsize=0.18cm 2.0,arrowlength=1.4,%
arrowinset=0.4]{<-}(6.7171874,-3.114165)(7.6371875,-0.8741651)
\psline[linewidth=0.04cm,linecolor=color98b,arrowsize=0.18cm 2.0,%
arrowlength=1.4,arrowinset=0.4]{->}(6.6971,-3.154)(7.8971,1.7858)
\psline[linewidth=0.04cm,linecolor=color97b,arrowsize=0.18cm 2.0,%
arrowlength=1.4,arrowinset=0.4]{->}(7.8971,1.7858)(8.657187,-0.5141651)
\psline[linewidth=0.04cm,arrowsize=0.18cm 2.0,%
arrowlength=1.4,arrowinset=0.4]{->}(8.657,-0.53416)(8.077188,-2.9541652)
\psline[linewidth=0.04cm,linecolor=color98b,arrowsize=0.18cm 2.0,%
arrowlength=1.4,arrowinset=0.4]{->}(8.057,-2.9941)(8.817,1.5858)
\psline[linewidth=0.04cm,linecolor=color97b,arrowsize=0.18cm 2.0,%
arrowlength=1.4,arrowinset=0.4]{->}(8.817,1.6058)(9.377,-0.31416)
\psline[linewidth=0.04cm,arrowsize=0.18cm 2.0,arrowlength=1.4,%
arrowinset=0.4]{->}(9.377,-0.33416)(9.02,-2.41)
\psline[linewidth=0.04cm,linecolor=color98b,arrowsize=0.18cm 2.0,%
arrowlength=1.4,arrowinset=0.4]{->}(9.0,-2.47)(9.497,1.5058)
\rput(4.06,-2.42){\huge $y_1$}
\rput(6.44,-3.54){\huge $y_2$}
\rput(7.94,-3.40){\huge $y_3$}
\rput(9.76,-0.78){\huge $x_4$}
\rput(9.35,-2.84){\huge $y_4$}
\rput(10.8,-2.0){\huge $C$}
\rput(14.3,3.4){\huge $D$}
\rput(0.19,-0.20){\huge $y_0$}
\rput(12.0,-0.3){\huge $x_\infty$}
\rput(2.76,1.99){\huge $x_0$}
\rput(6.65,-0.60){\huge $x_2$}
\rput(7.99,-1.11){\huge $x_3$}
\rput(5.19,0.28){\huge $x_1$}
\psdots[dotsize=0.18](8.06,-2.99)
\psdots[dotsize=0.18](5.58,-0.014)
\psdots[dotsize=0.18](6.70,-3.13)
\psdots[dotsize=0.18](2.98,1.59)
\psdots[dotsize=0.18](0.24,0.22)
\psdots[dotsize=0.18](4.60,-2.33)
\psdots[dotsize=0.18](0.24,0.20)
\psdots[dotsize=0.18](9.23,-0.57)
\psdots[dotsize=0.18](8.42,-0.75)
\psdots[dotsize=0.18](7.30,-0.65)
\psdots[dotsize=0.18](9.00,-2.47)
\psdots[dotsize=0.09](10.4,-0.50)
\psdots[dotsize=0.09](10.6,-0.45)
\psdots[dotsize=0.09](10.8,-0.40)
\end{pspicture} 
}

\caption{Forward-backward versus Douglas-Rachford: 
As in Example~\ref{plcex:2}, let $C$ and $D$ be two closed convex sets 
and consider the problem \eqref{plce:11} of finding a point $x_\infty$
in $C$ at minimum distance from $D$. Let us set $f_1=\iota_C$ and 
$f_2=d_D^2/2$. 
Top: The forward-backward algorithm 
with $\gamma_n\equiv 1$ and $\lambda_n\equiv 1$.
As seen in Example~\ref{plcex:2}, it assumes the form of the 
alternating projection method \eqref{plce:14}. 
Bottom: The Douglas-Rachford algorithm with $\gamma=1$ and
$\lambda_n\equiv 1$.
Table~\ref{plct:prop}.\ref{plct:propviii} yields 
${\mathrm{prox}}_{f_1}=P_C$ and
Table~\ref{plct:prop}.\ref{plct:propA380} yields 
${\mathrm{prox}}_{f_2}\colon x \mapsto (x+P_Dx)/2$. 
Therefore the updating rule in
Algorithm~\ref{plcalgo:DR} reduces to $x_n=(y_n+P_Dy_n)/2$ and 
$y_{n+1}=P_C(2x_n-y_n)+y_n-x_n=P_C(P_Dy_n)+y_n-x_n$.
}
\label{plcfig:A380}
\end{center}
\end{figure}

Applications of the Douglas-Rachford algorithm to signal and image
processing can be found in
\cite{plcChau09,plcJsts07,plcDupe09,plcDura08,%
plcSetz09,plcSetz10,plcStei09}.

The limiting case of the Douglas-Rachford
algorithm in which $\lambda_n\equiv 2$ is the
\emph{Peaceman-Rachford algorithm} \cite{plcOpti04,plcEcks92,plcLion79}. 
Its convergence requires additional assumptions (for instance, that
$f_2$ be strictly convex and real-valued) \cite{plcJoca09}.

\section{Dykstra-like splitting}
\label{plcsec:5}

In this section we consider problems involving a quadratic term
penalizing the deviation from a reference signal $r$.

\begin{problem}
\label{plcprob:4}
Let $f$ and $g$ be functions in $\Gamma_0({\mathbb R}^N)$ such that
${\mathrm{dom}}\, f\cap{\mathrm{dom}}\, g\neq\varnothing$, 
and let $r\in{\mathbb R}^N$. 
The problem is to 
\begin{equation}
\label{plce:222}
\underset{x\in{\mathbb R}^N}{\text{minimize}}\;\;
f(x)+g(x)+\frac12\|x-r\|^2.
\end{equation}
\end{problem}

It follows at once from \eqref{plce:prox1} that Problem~\ref{plcprob:4} 
admits a unique solution, namely $x={\mathrm{prox}}_{f+g}\,r$. 
Unfortunately, the proximity operator of the sum of two
functions is usually intractable. To compute it iteratively, 
we can observe that \eqref{plce:222} can be viewed as an instance of 
\eqref{plce:3} in Problem~\ref{plcprob:3} with $f_1=f$ and 
$f_2=g+\|\cdot-r\|^2/2$. However, in this Douglas-Rachford framework, the 
additional qualification condition \eqref{plce:cq1} needs to be imposed.
In the present setting we require only the minimal feasibility condition 
${\mathrm{dom}}\, f\cap{\mathrm{dom}}\, g\neq\varnothing$. 

\begin{algorithm}[Dykstra-like proximal algorithm]
\label{plcalgo:D}
~\\
Set $x_0=r$, $p_0=0$, $q_0=0$\\
For\;$n=0,1,\ldots$\\[-3mm]
\begin{equation}
\hskip -8.3cm
\begin{array}{l}
\left\lfloor
\begin{array}{l}
y_n={\mathrm{prox}}_g(x_n+p_n)\\
p_{n+1}=x_n+p_n-y_n\\
x_{n+1}={\mathrm{prox}}_f(y_n+q_n)\\
q_{n+1}=y_n+q_n-x_{n+1}.
\end{array}
\right.
\end{array}
\end{equation}
\end{algorithm}

\begin{proposition} {\rm\cite{plcPjo208}}
\label{plcp:6}
Every sequence $(x_n)_{n\in{\mathbb N}}$ generated by 
Algorithm~\ref{plcalgo:D} converges to the solution to 
Problem~\ref{plcprob:4}.
\end{proposition}

\begin{plcexample}[best approximation]
\label{plcex:10}
Let $f$ and $g$ be the indicator functions of closed convex
sets $C$ and $D$, respectively, in Problem~\ref{plcprob:4}. Then the
problem is to find the best approximation to $r$ from $C\cap D$, i.e., 
the projection of $r$ onto $C\cap D$. In this case, 
since ${\mathrm{prox}}_f=P_C$ and ${\mathrm{prox}}_g=P_D$, the above
algorithm reduces to Dykstra's projection method 
\cite{plcBoyl86,plcDyks83}.
\end{plcexample}

\begin{plcexample}[denoising]
\label{plcex:11}
Consider the problem of recovering a signal $\overline{x}$ from 
a noisy observation $r=\overline{x}+w$, where $w$ models noise. 
If $f$ and $g$ are functions in $\Gamma_0({\mathbb R}^N)$ promoting 
certain properties of $\overline{x}$, adopting a least-squares data 
fitting objective leads to the variational denoising problem 
\eqref{plce:222}.
\end{plcexample}

\section{Composite problems}
\label{plcsec:6}

We focus on variational problems with $m=2$ functions
involving explicitly a linear transformation.

\begin{problem}
\label{plcprob:7}
Let $f\in\Gamma_0({\mathbb R}^{N})$, let 
$g\in\Gamma_0({\mathbb R}^{M})$, and 
let $L\in{\mathbb R}^{M\times N}\smallsetminus\{0\}$ be such that 
${\mathrm{dom}}\, g\cap L({\mathrm{dom}}\, f)\neq\varnothing$ 
and $f(x)+g(Lx)\to+\infty$ as $\|x\|\to+\infty$. The problem is to
\begin{equation}
\label{plce:7}
\underset{x\in{\mathbb R}^N}{\text{minimize}}\;\;f(x)+g(Lx).
\end{equation}
\end{problem}

Our assumptions guarantee that Problem~\ref{plcprob:7} possesses at least
one solution. To find such a solution, several scenarios can be
contemplated.

\subsection{Forward-backward splitting}

Suppose that in Problem~\ref{plcprob:7} 
$g$ is differentiable with a $\tau$-Lipschitz continuous
gradient (see \eqref{plce:4}).
Now set $f_1=f$ and $f_2=g\circ L$. 
Then $f_2$ is differentiable and its gradient 
\begin{equation}
\nabla f_2=L^\top\circ\nabla g\circ L
\end{equation}
is $\beta$-Lipschitz continuous, with $\beta=\tau\|L\|^2$.
Hence, we can apply the forward-backward splitting method,
as implemented in Algorithm~\ref{plcalgo:FB}. As seen in 
\eqref{plce:FB}, it 
operates with the updating rule
\begin{equation}
\begin{array}{l}
\left\lfloor
\begin{array}{l}
\gamma_n\in [\varepsilon,2/(\tau\|L\|^2)-\varepsilon]\\
y_n=x_n-\gamma_nL^\top\nabla g(Lx_n)\\
\lambda_n \in [\varepsilon,1]\\
x_{n+1}=x_n+\lambda_n({\mathrm{prox}}_{\gamma_n f}y_n-x_n).
\end{array}
\right.\\
\end{array}
\end{equation}
Convergence is guaranteed by Proposition~\ref{plcp:2}.

\subsection{Douglas-Rachford splitting}

Suppose that in Problem~\ref{plcprob:7} the matrix $L$ satisfies
\begin{equation}
\label{plce:rio2009-08-30a}
LL^\top=\nu I,\quad\text{where}\quad\nu\in\,\left]0,+\infty\right[
\end{equation}
and $({\mathrm{ri}}\,{\mathrm{dom}}\, g)\cap
{\mathrm{ri}}\, L({\mathrm{dom}}\, f)
\neq\varnothing$.
Let us set $f_1=f$ and $f_2=g\circ L$. As seen in 
Table~\ref{plct:prop}.\ref{plct:prop8}, ${\mathrm{prox}}_{f_2}$ 
has a closed-form expression in terms of ${\mathrm{prox}}_{g}$ 
and we can therefore apply the Douglas-Rachford splitting method
(Algorithm~\ref{plcalgo:DR}). In this scenario, the updating 
rule reads
\begin{equation}
\begin{array}{l}
\left\lfloor
\begin{array}{l}
x_n=y_n+\nu^{-1}L^\top
\big({\mathrm{prox}}_{\gamma\nu g}(Ly_n)-Ly_n\big)\\
\lambda_n\in [\varepsilon,2-\varepsilon]\\
y_{n+1}=y_n+\lambda_n
\big({\mathrm{prox}}_{\gamma f}\big(2x_n-y_n\big)-x_n\big). 
\end{array}
\right.
\end{array}
\end{equation}
Convergence is guaranteed by Proposition~\ref{plcp:3}.

\subsection{Dual forward-backward splitting}

Suppose that in Problem~\ref{plcprob:7} $f=h+\|\cdot-r\|^2/2$, where 
$h\in\Gamma_0({\mathbb R}^N)$ and $r\in{\mathbb R}^N$. 
Then \eqref{plce:7} becomes
\begin{equation}
\label{plce:77}
\underset{x\in{\mathbb R}^N}{\text{minimize}}\;\;h(x)+g(Lx)
+\frac12\|x-r\|^2,
\end{equation}
which models various signal recovery problems, e.g.,
\cite{plcCham04,plcCham05,plcViet09,plcDida09,plcPott93,plcYoul78}. 
If \eqref{plce:rio2009-08-30a} holds, 
${\mathrm{prox}}_{g\circ L}$ is decomposable, 
and \eqref{plce:77} can be solved
with the Dykstra-like method of Section~\ref{plcsec:5}, where
$f_1=h+\|\cdot-r\|^2/2$ (see Table~\ref{plct:prop}.\ref{plct:propLeblon})
and $f_2=g\circ L$ (see Table~\ref{plct:prop}.\ref{plct:prop8}).
Otherwise, we can exploit the nice properties of the 
Fenchel-Moreau-Rockafellar dual 
of \eqref{plce:77}, solve this dual problem by forward-backward splitting,
and recover the unique solution to \eqref{plce:77} \cite{plcViet09}. 

\begin{algorithm}[Dual forward-backward algorithm]
\label{plce:DFB}
~\\
Fix $\varepsilon\in\left]0,\min\{1,1/\|L\|^2\}\right[$,
$u_0\in{\mathbb R}^M$\\
For\;$n=0,1,\ldots$
\begin{equation}
\begin{array}{l}
\hskip -5.3cm
\left\lfloor
\begin{array}{l}
x_n={\mathrm{prox}}_h(r-L^\top u_n)\\[1mm]
\gamma_n\in\left[\varepsilon,2/\|L\|^{2}-\varepsilon\right]\\[1mm]
\lambda_n\in\left[\varepsilon,1\right]\\
u_{n+1}=u_n+\lambda_n\big({\mathrm{prox}}_{\gamma_n g^*}(u_n
+\gamma_nLx_n)-u_n\big).
\end{array}
\right.\\[2mm]
\end{array}
\end{equation}
\end{algorithm}

\begin{proposition}{\rm \cite{plcViet09}}
\label{plcp:9}
Assume that $({\mathrm{ri}}\,{\mathrm{dom}}\, g)\cap{\mathrm{ri}}\, 
L({\mathrm{dom}}\, h)\neq\varnothing$. Then every 
sequence $(x_n)_{n\in{\mathbb N}}$ generated by the dual forward-backward
algorithm \ref{plce:DFB} converges to the solution to \eqref{plce:77}.
\end{proposition}

\subsection{Alternating-direction method of multipliers}
\label{plcsec:6d}
Augmented Lagrangian techniques are classical approaches for solving 
Problem~\ref{plcprob:7} \cite{plcGaba76,plcGlow75} 
(see also \cite{plcGlow83,plcGlow89}).
First, observe that \eqref{plce:7} is equivalent to 
\begin{equation}
\label{plce:rio2009-08-30b}
\underset{\substack{x\in{\mathbb R}^N,\:y\in{\mathbb R}^M\\Lx=y}}
{\text{minimize}}\;\;f(x)+g(y).
\end{equation}
The \emph{augmented Lagrangian} of index 
$\gamma\in\,\left]0,+\infty\right[$ 
associated with \eqref{plce:rio2009-08-30b} is the saddle function
\begin{align}
\label{plce:rio2009-08-30c}
{\mathcal L}_\gamma\colon{\mathbb R}^N\times
{\mathbb R}^M\times{\mathbb R}^M&\to\,\left]-\infty,+\infty\right]
\nonumber\\
(x,y,z)&\mapsto f(x)+g(y)+
\frac{1}{\gamma}{z}^\top{(Lx-y)}+\frac{1}{2\gamma}\|Lx-y\|^2.
\end{align}
The alternating-direction method of multipliers consists in minimizing
${\mathcal L}_\gamma$ over $x$, then over $y$, and then applying a
proximal maximization step with respect to the Lagrange multiplier $z$.
Now suppose that 
\begin{equation}
\label{plce:rio2009-08-30d}
L^\top L\;\:\text{is invertible and}\;\:
({\mathrm{ri}}\,{\mathrm{dom}}\, g)\cap{\mathrm{ri}}\, 
L({\mathrm{dom}}\, f)\neq\varnothing.
\end{equation}
By analogy with \eqref{plce:prox1}, if we denote by 
${\mathrm{prox}}_f^L$ the 
operator which maps a point $y\in{\mathbb R}^M$ to the unique 
minimizer of $x\mapsto f(x)+\|Lx-y\|^2/2$, we obtain the following 
implementation. 

\begin{algorithm}[Alternating-direction method of multipliers (ADMM)]
\label{plce:ADMM}
~\\
Fix $\gamma>0$, $y_0\in{\mathbb R}^M$, $z_0\in{\mathbb R}^M$\\
For\;$n=0,1,\ldots$
\begin{equation}
\hskip -8.1cm
\begin{array}{l}
\left\lfloor
\begin{array}{l}
x_{n}={\mathrm{prox}}_{\gamma f}^L(y_n-z_n)\\
s_n=Lx_n\\
y_{n+1}={\mathrm{prox}}_{\gamma g}(s_n+z_n)\\
z_{n+1}=z_n+s_n-y_{n+1}.
\end{array}
\right.\\[2mm]
\end{array}
\end{equation}
\end{algorithm}

The convergence of the sequence $(x_n)_{n\in{\mathbb N}}$ thus produced
under assumption \eqref{plce:rio2009-08-30d} has been investigated 
in several places, e.g., \cite{plcGlow83,plcGaba76,plcGlow89}. 
It was first observed 
in \cite{plcGaba83} that the ADMM algorithm can be derived
from an application of the Douglas-Rachford algorithm to the
dual of \eqref{plce:7}. This analysis was pursued in \cite{plcEcks92},
where the convergence of $(x_n)_{n\in{\mathbb N}}$ to a solution to
\eqref{plce:7} is shown. Variants of the method relaxing the 
requirements on $L$ in \eqref{plce:rio2009-08-30d} have been 
proposed \cite{plcAtto09,plcChen94}.

In image processing, ADMM 
was applied in \cite{plcGold09} to an $\ell_1$ regularization 
problem under the name ``alternating split Bregman algorithm.'' 
Further applications and connections are found in 
\cite{plcFigu09b,plcEsse09,plcSetz09,plcZhan09}.

\section{Problems with $m\geq 2$ functions}
\label{plcsec:7}

We return to the general minimization problem \eqref{plce:1}.

\begin{problem}
\label{plcprob:24}
Let $f_1$,\ldots,$f_m$ be functions in $\Gamma_0({\mathbb R}^N)$ 
such that
\begin{equation}
\label{plce:rio2009-09-01a}
({\mathrm{ri}}\,{\mathrm{dom}}\, f_1)\cap\cdots
\cap({\mathrm{ri}}\,{\mathrm{dom}}\, f_m)\neq\varnothing
\end{equation}
and $f_1(x)+\cdots+f_m(x)\to+\infty$ as $\|x\|\to+\infty$. 
The problem is to
\begin{equation}
\label{plce:rio2009-08-29a}
\underset{x\in{\mathbb R}^N}{\text{minimize}}\;\;f_1(x)+\cdots+f_m(x).
\end{equation}
\end{problem}

Since the methods described so far are designed for $m=2$ functions,
we can attempt to reformulate \eqref{plce:rio2009-08-29a} as a 
2-function problem in the $m$-fold product space
\begin{equation}
\label{plce:rio2009-09-01b}
{\boldsymbol{\mathcal H}}={\mathbb R}^N\times\cdots\times{\mathbb R}^N 
\end{equation}
(such techniques
were introduced in \cite{plcPier76,plcPier84} and have been 
used in the context of convex feasibility
problems in \cite{plcBaus96,plcSign94,plcImag97}). 
To this end, observe that
\eqref{plce:rio2009-08-29a} can be rewritten in 
${\boldsymbol{\mathcal H}}$ as
\begin{equation}
\label{plce:rio2009-08-29b}
\underset{\substack{(x_1,\ldots,x_m)\in{\boldsymbol{\mathcal H}}\\
x_1=\cdots=x_m}}{\text{minimize}}\;\;f_1(x_1)+\cdots+f_m(x_m).
\end{equation}
If we denote by ${\boldsymbol x}=(x_1,\ldots,x_m)$ a 
generic element in ${\boldsymbol{\mathcal H}}$, 
\eqref{plce:rio2009-08-29b} 
is equivalent to 
\begin{equation}
\label{plce:rio2009-08-29c}
\underset{{\boldsymbol x}\in{\boldsymbol{\mathcal H}}}
{\text{minimize}}\;\;\iota_{\boldsymbol D}({\boldsymbol x})
+{\boldsymbol f}({\boldsymbol x}),
\end{equation}
where 
\begin{equation}
\begin{cases}
{\boldsymbol D}=\big\{(x,\ldots,x)\in{\boldsymbol{\mathcal H}}\:|\:
x\in{\mathbb R}^N\big\}\\
{\boldsymbol f}\colon{\boldsymbol x}\mapsto f_1(x_1)+\cdots+f_m(x_m).
\end{cases}
\end{equation}
We are thus back to a problem involving two functions in the larger
space ${\boldsymbol{\mathcal H}}$. 
In some cases, this observation makes it possible to
obtain convergent methods from the algorithms discussed in the preceding
sections. For instance, the following parallel algorithm was 
derived from the Douglas-Rachford algorithm in \cite{plcInvp08}
(see also \cite{plcJoca09} for further analysis and connections 
with Spingarn's splitting method \cite{plcSpin83}).

\begin{algorithm}[Parallel proximal algorithm (PPXA)]
\label{plcalgo:ppxa}
~\\
Fix $\varepsilon\in\:]0,1[,\;\gamma>0,$
$(\omega_{i})_{1\leq i\leq m}\in\left]0,1\right]^m$ such that\\[1mm]
$\sum_{i=1}^m\omega_i=1$,\;
$y_{1,0}\in{\mathbb R}^N,\ldots,y_{m,0}\in{\mathbb R}^N$\\[1mm]
Set $x_0={\sum_{i=1}^m}\,\omega_iy_{i,0}$

\hskip -4.5mm
$
\begin{array}{l}
\text{For}\;n=0,1,\ldots\\
\left\lfloor
\begin{array}{l}
\text{For}\;i=1,\ldots,m\\
\quad\left\lfloor
\begin{array}{l}
p_{i,n}={\mathrm{prox}}_{\gamma f_i/\omega_i}y_{i,n}\\
\end{array}
\right.\\[1mm]
p_n=\displaystyle{\sum_{i=1}^m}\,\omega_ip_{i,n}\\
\varepsilon\leq\lambda_n\leq 2-\varepsilon\\
\text{For}\;i=1,\ldots,m\\
\quad\left\lfloor
\begin{array}{l}
y_{i,n+1}=y_{i,n}+\lambda_n\big(2p_n-x_n-p_{i,n}\big)
\end{array}
\right.\\[2mm]
x_{n+1}=x_n+\lambda_n(p_{n}-x_n).
\end{array}
\right.\\
\end{array}
$
\end{algorithm}

\begin{proposition}{\rm \cite{plcInvp08}}
\label{plcp:19}
Every sequence $(x_n)_{n\in{\mathbb N}}$ generated by
Algorithm~\ref{plcalgo:ppxa} converges to a solution to 
Problem~\ref{plcprob:24}.
\end{proposition}

\begin{plcexample}[image recovery]
\label{plcex:33333}
In many imaging problems, we record an observation 
$y\in {\mathbb R}^M$ of an image $\overline{z}\in {\mathbb R}^K$ 
degraded by a matrix $L\in {\mathbb R}^{M\times K}$ and 
corrupted by noise.
In the spirit of a number of recent investigations 
(see \cite{plcChau07} and the references therein), 
a tight frame representation of the 
images under consideration can be used. 
This representation is defined through a synthesis matrix 
$F^\top\in{\mathbb R}^{K\times N}$ (with $K \leq N$) such that 
$F^\top F=\nu I$, for some $\nu\in\,\left]0,+\infty\right[$. 
Thus, the original image can be written as 
$\overline{z}=F^\top\overline{x}$, where 
$\overline{x}\in {\mathbb R}^N$ is a vector of
frame coefficients to be estimated. For this purpose,
we consider the problem 
\begin{equation}\label{plce:explus}
\underset{x\in C}{\text{minimize}}\;\;
\frac12\|LF^\top x-y\|^2+\Phi(x)+\operatorname{tv}(F^\top x),
\end{equation}
where $C$ is a closed convex set modeling a constraint on
$\overline{z}$, the quadratic term is the standard least-squares
data fidelity term, $\Phi$ is a real-valued convex function on 
${\mathbb R}^N$ (e.g., a weighted $\ell^1$ norm) 
introducing a regularization on the frame coefficients, and 
$\operatorname{tv}$ is a discrete total variation function 
aiming at preserving piecewise smooth areas and sharp 
edges \cite{plcRudi92}.
Using appropriate gradient filters in the computation of 
$\operatorname{tv}$, it is possible to decompose it as 
a sum of convex functions $(\operatorname{tv}_i)_{1\leq i\leq q}$, 
the proximity operators
of which can be expressed in closed form \cite{plcInvp08,plcPust09}.
Thus, \eqref{plce:explus} appears as a special case of  
\eqref{plce:rio2009-08-29a}
with $m = q+3$, $f_1=\iota_C$, $f_2=\|LF^\top \cdot-y\|^2/2$,
$f_3=\Phi$, and $f_{3+i}=\operatorname{tv}_i(F^\top \cdot)$ for
$i\in\{1,\ldots,q\}$. Since a tight frame is employed,
the proximity operators of $f_2$ and $(f_{3+i})_{1\leq i \leq q}$
can be deduced from Table~\ref{plct:prop}.\ref{plct:prop8}.
Thus, the PPXA algorithm is well suited for solving
this problem numerically.
\end{plcexample}

A product space strategy can also be adopted to 
address the following extension of Problem~\ref{plcprob:4}.
\begin{problem}
\label{plcprob:52}
Let $f_1$, \ldots, $f_m$ be functions in $\Gamma_0({\mathbb R}^N)$ 
such that
${\mathrm{dom}}\, f_1\cap\cdots\cap{\mathrm{dom}}\, f_m\neq\varnothing$,
let $(\omega_i)_{1\leq i \leq m}\in\:]0,1]^m$ be such 
that $\sum_{i=1}^m\omega_i=1$, and let $r\in{\mathbb R}^N$. 
The problem is to
\begin{equation}
\label{plce:manille-16mai2008-1}
\underset{x\in{\mathbb R}^N}{\mathrm{minimize}}\;\;
\sum_{i=1}^m\omega_if_i(x)
+\frac12\|x-r\|^2.
\end{equation}
\end{problem}

\begin{algorithm}[Parallel Dykstra-like proximal algorithm]
\label{plcalgo:PDL}
~\\
Set $x_0=r$, $z_{1,0}=x_0$, \ldots, $z_{m,0}=x_0$\\
For\;$n=0,1,\ldots$
\begin{equation}
\hskip -7.5cm
\begin{array}{l}
\left\lfloor
\begin{array}{l}
\text{For}\;i=1,\ldots,m\\
\left\lfloor
\begin{array}{l}
p_{i,n}={\mathrm{prox}}_{f_i}z_{i,n}
\end{array}
\right.\\[1mm]
x_{n+1}={\sum_{i=1}^m}\,\omega_i p_{i,n}\\
\text{For}\;i=1,\ldots,m\\
\left\lfloor
\begin{array}{l}
z_{i,n+1}=x_{n+1}+z_{i,n}-p_{i,n}.
\end{array}
\right.\\[2mm]
\end{array}
\right.
\end{array}
\end{equation}
\end{algorithm}

\begin{proposition}{\rm \cite{plcJoca09}}
\label{plcp:15}
Every sequence $(x_n)_{n\in{\mathbb N}}$ generated by 
Algorithm~\ref{plcalgo:PDL} converges to 
the solution to Problem~\ref{plcprob:52}.
\end{proposition}

Next, we consider a composite problem.

\begin{problem}
\label{plcprob:53}
For every $i\in\{1,\ldots,m\}$, let 
$g_i\in\Gamma_0({\mathbb R}^{M_i})$ and 
let $L_i\in{\mathbb R}^{M_i\times N}$. Assume that
\begin{equation}
(\exists\, q\in{\mathbb R}^N)\quad L_1q\in{\mathrm{ri}}\,{\mathrm{dom}}\, 
g_1,\ldots, L_mq\in{\mathrm{ri}}\,{\mathrm{dom}}\, g_m,
\end{equation}
that $g_1(L_1x)+\cdots+g_m(L_mx)\to +\infty$ as $\|x\|\to+\infty$, and
that $Q=\sum_{1\leq i \leq m} L_i^\top L_i$ is invertible.
The problem is to 
\begin{equation}
\label{plce:32}
\underset{x\in{\mathbb R}^N}{\text{minimize}}
\;\;g_1(L_1x)+\cdots+g_m(L_mx).
\end{equation}
\end{problem}

Proceeding as in \eqref{plce:rio2009-08-29b} and 
\eqref{plce:rio2009-08-29c}, \eqref{plce:32} can be recast as
\begin{equation}
\underset{\substack{{\boldsymbol x}\in{\boldsymbol{\mathcal H}},
\:{\boldsymbol y}\in \boldsymbol{\mathcal G}\\
{\boldsymbol y} = {\boldsymbol L}{\boldsymbol x}}}
{\text{minimize}}\;\;\iota_{\boldsymbol D}({\boldsymbol x})
+{\boldsymbol g}({\boldsymbol y}),
\end{equation}
where 
\begin{equation}
\begin{cases}
{\boldsymbol{\mathcal H}}={\mathbb R}^{N}
\times\cdots\times{\mathbb R}^{N},\:
\boldsymbol{\mathcal G}={\mathbb R}^{M_1}\times\cdots\times
{\mathbb R}^{M_m}\\
{\boldsymbol L}\colon{\boldsymbol{\mathcal H}}\to
\boldsymbol{\mathcal G} \colon {\boldsymbol x}
\mapsto (L_1x_1,\ldots,L_mx_m)\\
{\boldsymbol g}\colon\boldsymbol{\mathcal G}\to\,
\left]-\infty,+\infty\right]\colon{\boldsymbol y}
\mapsto g_1(y_1)+\cdots+g_m(y_m).
\end{cases}
\end{equation}
In turn, a solution to \eqref{plce:32} can be obtained as the limit of
the sequence $(x_n)_{n\in{\mathbb N}}$ constructed by the 
following algorithm, which can be derived from the alternating-direction 
method of multipliers of Section~\ref{plcsec:6d} (alternative parallel 
offsprings of ADMM exist, see for instance \cite{plcEcks94}).

\begin{algorithm}[Simultaneous-direction method of multipliers (SDMM)]
\label{plcalgo:vmm}
~\\
Fix $\gamma>0$,
$\:y_{1,0}\in{\mathbb R}^{M_1},\ldots,\,y_{m,0}\in{\mathbb R}^{M_m}$, 
$z_{1,0}\in{\mathbb R}^{M_1},\ldots,\,z_{m,0}\in{\mathbb R}^{M_m}$
\begin{equation}
\label{plce:vmm}
\hskip -6.8cm
\begin{array}{l}
\text{For}\;n=0,1,\ldots\\
\left\lfloor
\begin{array}{l}
x_{n}=Q^{-1}\sum_{i=1}^m L_i^\top (y_{i,n}-z_{i,n})\\
\mathrm{For}\;i=1,\ldots,m\\
\quad\left\lfloor
\begin{array}{l}
s_{i,n} = L_i x_{n}\\
y_{i,n+1}={\mathrm{prox}}_{\gamma g_i}(s_{i,n}+z_{i,n})\\
z_{i,n+1} = z_{i,n}+s_{i,n}-y_{i,n+1}\\[1mm]
\end{array}
\right.
\end{array}
\right.
\end{array}
\end{equation}
\end{algorithm}

This algorithm was derived from a slightly different viewpoint in 
\cite{plcSetz10} with a connection with the work of \cite{plcFigu09}. 
In these papers, SDMM is applied to deblurring in the presence 
of Poisson noise. 
The computation of $x_n$ in \eqref{plce:vmm} requires the 
solution of a positive-definite symmetric system of linear 
equations. Efficient methods for solving such systems can be found in 
\cite{plcGolu96}. In certain situations, fast Fourier 
diagonalization is also an option \cite{plcFigu09b,plcFigu09}.

In the above algorithms, the proximal vectors, as well as the 
auxiliary vectors, can be computed simultaneously at each iteration. 
This parallel structure is useful when the algorithms are implemented on 
multicore architectures. A parallel proximal algorithm is also
available to solve multicomponent signal processing problems  
\cite{plcLuis09}. This framework captures in particular problem
formulations found in 
\cite{plcAujo04,plcAujo05,plcGold85,plcHuan08,plcVese04}.
Let us add that an alternative splitting
framework applicable to \eqref{plce:rio2009-08-29a}
was recently proposed in \cite{plcSvai09}.

\section{Conclusion}
\label{plcsec:9}
We have presented a panel of convex optimization 
algorithms sharing two main features. First, they employ proximity 
operators, a powerful generalization of the notion of a projection 
operator. Second, they operate by splitting the objective to
be minimized into simpler functions that are dealt with individually.
These methods are applicable to a wide class of signal and image 
processing problems ranging from restoration and reconstruction to 
synthesis and design.
One of the main advantages of these algorithms is that they can be used 
to minimize nondifferentiable objectives, such as those commonly
encountered in sparse approximation and compressed sensing, or
in hard-constrained problems. 
Finally, let us note that the variational problems described in
\eqref{plce:222}, \eqref{plce:77}, and 
\eqref{plce:manille-16mai2008-1},
consist of computing a proximity operator. Therefore the associated
algorithms can be used as a subroutine to compute approximately 
proximity operators within a proximal splitting algorithm, 
provided the latter is error tolerant (see 
\cite{plcOpti04,plcJoca09,plcViet09,plcEcks92,plcRock76} 
for convergence properties under approximate proximal computations). 
An application of this principle can be found in \cite{plcChau09}.


\begin{thebibliography}{100}
\providecommand{\url}[1]{{#1}}
\providecommand{\urlprefix}{URL }
\expandafter\ifx\csname urlstyle\endcsname\relax
  \providecommand{\doi}[1]{DOI~\discretionary{}{}{}#1}\else
  \providecommand{\doi}{DOI~\discretionary{}{}{}\begingroup
  \urlstyle{rm}\Url}\fi

\bibitem{plcAcke80}
Acker, F., Prestel, M.A.: Convergence d'un sch\'ema de minimisation altern\'ee.
\newblock Ann. Fac. Sci. Toulouse V. S\'er. Math. \textbf{2}, 1--9 (1980)

\bibitem{plcFigu09b}
Afonso, M.V., Bioucas-Dias, J.M., Figueiredo, M.A.T.: Fast image recovery using
  variable splitting and constrained optimization (2009).
\newblock {http://arxiv.org/abs/0910.4887}

\bibitem{plcAnto01}
Antoniadis, A., Fan, J.: Regularization of wavelet approximations.
\newblock J. Amer. Statist. Assoc. \textbf{96}, 939--967 (2001)

\bibitem{plcAnto02}
Antoniadis, A., Leporini, D., Pesquet, J.C.: Wavelet thresholding for some
  classes of non-{G}aussian noise.
\newblock Statist. Neerlandica \textbf{56}, 434--453 (2002)

\bibitem{plcAtto09}
Attouch, H., Soueycatt, M.: Augmented {L}agrangian and proximal alternating
  direction methods of multipliers in {H}ilbert spaces -- applications to
  games, {PDE}s and control.
\newblock Pacific J. Optim. \textbf{5}, 17--37 (2009)

\bibitem{plcAubi98}
Aubin, J.P.: Optima and Equilibria -- An Introduction to Nonlinear Analysis,
  2nd edn.
\newblock Springer-Verlag, New York (1998)

\bibitem{plcAujo04}
Aujol, J.F., Aubert, G., Blanc-F\'eraud, L., Chambolle, A.: Image decomposition
  into a bounded variation component and an oscillating component.
\newblock J. Math. Imaging Vision \textbf{22}, 71--88 (2005)

\bibitem{plcAujo05}
Aujol, J.F., Chambolle, A.: Dual norms and image decomposition models.
\newblock Int. J. Computer Vision \textbf{63}, 85--104 (2005)

\bibitem{plcBaus94}
Bauschke, H.H., Borwein, J.M.: Dykstra's alternating projection algorithm for
  two sets.
\newblock J. Approx. Theory \textbf{79}, 418--443 (1994)

\bibitem{plcBaus96}
Bauschke, H.H., Borwein, J.M.: On projection algorithms for solving convex
  feasibility problems.
\newblock SIAM Rev. \textbf{38}, 367--426 (1996)

\bibitem{plcMoor01}
Bauschke, H.H., Combettes, P.L.: A weak-to-strong convergence principle for
  {F}ej\'er-monotone methods in {H}ilbert spaces.
\newblock Math. Oper. Res. \textbf{26}, 248--264 (2001)

\bibitem{plcPjo208}
Bauschke, H.H., Combettes, P.L.: A {D}ykstra-like algorithm for two monotone
  operators.
\newblock Pacific J. Optim. \textbf{4}, 383--391 (2008)

\bibitem{plcLivr10}
Bauschke, H.H., Combettes, P.L.: Convex Analysis and Monotone Operator Theory
  in {H}ilbert Spaces.
\newblock Springer-Verlag (2011).
\newblock To appear

\bibitem{plcReic05}
Bauschke, H.H., Combettes, P.L., Reich, S.: The asymptotic behavior of the
  composition of two resolvents.
\newblock Nonlinear Anal. \textbf{60}, 283--301 (2005)

\bibitem{plcBec09b}
Beck, A., Teboulle, M.: Fast gradient-based algorithms for constrained total
  variation image denoising and deblurring problems.
\newblock IEEE Trans. Image Process. \textbf{18}, 2419--2434 (2009)

\bibitem{plcBec09a}
Beck, A., Teboulle, M.: A fast iterative shrinkage-thresholding algorithm for
  linear inverse problems.
\newblock SIAM J. Imaging Sci. \textbf{2}, 183--202 (2009)

\bibitem{plcBect04}
Bect, J., Blanc-F\'eraud, L., Aubert, G., Chambolle, A.: A $\ell^1$ unified
  variational framework for image restoration.
\newblock Lecture Notes in Comput. Sci. \textbf{3024}, 1--13 (2004)

\bibitem{plcBenv10}
Benvenuto, F., Zanella, R., Zanni, L., Bertero, M.: Nonnegative least-squares
  image deblurring: improved gradient projection approaches.
\newblock Inverse Problems \textbf{26}, 18 (2010).
\newblock Art. 025004

\bibitem{plcBert97}
Bertsekas, D.P., Tsitsiklis, J.N.: Parallel and Distributed Computation:
  Numerical Methods.
\newblock Athena Scientific, Belmont, MA (1997)

\bibitem{plcBiou07}
Bioucas-Dias, J.M., Figueiredo, M.A.T.: A new {TwIST}: Two-step iterative
  shrinkage/thresholding algorithms for image restoration.
\newblock IEEE Trans. Image Process. \textbf{16}, 2992--3004 (2007)

\bibitem{plcBoum93}
Bouman, C., Sauer, K.: A generalized {G}aussian image model for edge-preserving
  {MAP} estimation.
\newblock IEEE Trans. Image Process. \textbf{2}, 296--310 (1993)

\bibitem{plcBoyl86}
Boyle, J.P., Dykstra, R.L.: A method for finding projections onto the
  intersection of convex sets in {H}ilbert spaces.
\newblock Lecture Notes in Statist. \textbf{37}, 28--47 (1986)

\bibitem{plcBred09}
Bredies, K.: A forward-backward splitting algorithm for the minimization of
  non-smooth convex functionals in {B}anach space.
\newblock Inverse Problems \textbf{25}, 20 (2009).
\newblock Art. 015005

\bibitem{plcBred08}
Bredies, K., Lorenz, D.A.: Linear convergence of iterative soft-thresholding.
\newblock J. Fourier Anal. Appl. \textbf{14}, 813--837 (2008)

\bibitem{plcBreg65}
Br\`{e}gman, L.M.: The method of successive projection for finding a common
  point of convex sets.
\newblock Soviet Math. Dokl. \textbf{6}, 688--692 (1965)

\bibitem{plcBrez78}
Br\'ezis, H., Lions, P.L.: Produits infinis de r\'esolvantes.
\newblock Israel J. Math. \textbf{29}, 329--345 (1978)

\bibitem{plcLuis09}
{Brice\~{n}o-Arias}, L.M., Combettes, P.L.: Convex variational formulation with
  smooth coupling for multicomponent signal decomposition and recovery.
\newblock Numer. Math. Theory Methods Appl. \textbf{2}, 485--508 (2009)

\bibitem{plcCaij09}
Cai, J.F., Chan, R.H., Shen, L., Shen, Z.: Convergence analysis of tight
  framelet approach for missing data recovery.
\newblock Adv. Comput. Math. \textbf{31}, 87--113 (2009)

\bibitem{plcChan09}
Cai, J.F., Chan, R.H., Shen, L., Shen, Z.: Simultaneously inpainting in image
  and transformed domains.
\newblock Numer. Math. \textbf{112}, 509--533 (2009)

\bibitem{plcCaij08}
Cai, J.F., Chan, R.H., Shen, Z.: A framelet-based image inpainting algorithm.
\newblock Appl. Comput. Harm. Anal. \textbf{24}, 131--149 (2008)

\bibitem{plcCens97}
Censor, Y., Zenios, S.A.: Parallel Optimization: Theory, Algorithms and
  Applications.
\newblock Oxford University Press, New York (1997)

\bibitem{plcChaa09}
Cha\^ari, L., Pesquet, J.C., Ciuciu, P., Benazza-Benyahia, A.: An iterative
  method for parallel {MRI} {SENSE}-based reconstruction in the wavelet domain
  (2009).
\newblock
  {http://www-syscom.univ-mlv.fr/$\sim$chaari/downloads/MEDIA\_Chaari.pdf}

\bibitem{plcCham04}
Chambolle, A.: An algorithm for total variation minimization and applications.
\newblock J. Math. Imaging Vision \textbf{20}, 89--97 (2004)

\bibitem{plcCham05}
Chambolle, A.: Total variation minimization and a class of binary {MRF} model.
\newblock Lecture Notes in Comput. Sci. \textbf{3757}, 136--152 (2005)

\bibitem{plcCham98}
Chambolle, A., DeVore, R.A., Lee, N.Y., Lucier, B.J.: Nonlinear wavelet image
  processing: Variational problems, compression, and noise removal through
  wavelet shrinkage.
\newblock IEEE Trans. Image Process. \textbf{7}, 319--335 (1998)

\bibitem{plcRaym08}
Chan, R.H., Setzer, S., Steidl, G.: Inpainting by flexible {H}aar-wavelet
  shrinkage.
\newblock SIAM J. Imaging Sci. \textbf{1}, 273--293 (2008)

\bibitem{plcChau07}
Chaux, C., Combettes, P.L., Pesquet, J.C., Wajs, V.R.: A variational
  formulation for frame-based inverse problems.
\newblock Inverse Problems \textbf{23}, 1495--1518 (2007)

\bibitem{plcChau09}
Chaux, C., Pesquet, J.C., Pustelnik, N.: Nested iterative algorithms for convex
  constrained image recovery problems.
\newblock SIAM J. Imaging Sci. \textbf{2}, 730--762 (2009)

\bibitem{plcChen94}
Chen, G., Teboulle, M.: A proximal-based decomposition method for convex
  minimization problems.
\newblock Math. Programming \textbf{64}, 81--101 (1994)

\bibitem{plcChen97}
Chen, G.H.G., Rockafellar, R.T.: Convergence rates in forward-backward
  splitting.
\newblock SIAM J. Optim. \textbf{7}, 421--444 (1997)

\bibitem{plcChen59}
Cheney, W., Goldstein, A.A.: Proximity maps for convex sets.
\newblock Proc. Amer. Math. Soc. \textbf{10}, 448--450 (1959)

\bibitem{plcProc93}
Combettes, P.L.: The foundations of set theoretic estimation.
\newblock Proc. IEEE \textbf{81}, 182--208 (1993)

\bibitem{plcSign94}
Combettes, P.L.: Inconsistent signal feasibility problems: Least-squares
  solutions in a product space.
\newblock IEEE Trans. Signal Process. \textbf{42}, 2955--2966 (1994)

\bibitem{plcAiep96}
Combettes, P.L.: The convex feasibility problem in image recovery.
\newblock In: P.~Hawkes (ed.) Advances in Imaging and Electron Physics,
  vol.~95, pp. 155--270. Academic Press, New York (1996)

\bibitem{plcImag97}
Combettes, P.L.: Convex set theoretic image recovery by extrapolated iterations
  of parallel subgradient projections.
\newblock IEEE Trans. Image Process. \textbf{6}, 493--506 (1997)

\bibitem{plcSmai01}
Combettes, P.L.: Convexit\'e et signal.
\newblock In: Proc. Congr\`es de Math\'ematiques Appliqu\'ees et Industrielles
  SMAI'01, pp. 6--16. Pompadour, France (2001)

\bibitem{plcSign03}
Combettes, P.L.: A block-iterative surrogate constraint splitting method for
  quadratic signal recovery.
\newblock IEEE Trans. Signal Process. \textbf{51}, 1771--1782 (2003)

\bibitem{plcOpti04}
Combettes, P.L.: Solving monotone inclusions via compositions of nonexpansive
  averaged operators.
\newblock Optimization \textbf{53}, 475--504 (2004)

\bibitem{plcJoca09}
Combettes, P.L.: Iterative construction of the resolvent of a sum of maximal
  monotone operators.
\newblock J. Convex Anal. \textbf{16}, 727--748 (2009)

\bibitem{plcSign99}
Combettes, P.L., Bondon, P.: Hard-constrained inconsistent signal feasibility
  problems.
\newblock IEEE Trans. Signal Process. \textbf{47}, 2460--2468 (1999)

\bibitem{plcViet09}
Combettes, P.L., {D\~ung}, D., {V\~u}, B.C.: Dualization of signal recovery
  problems (2009).
\newblock {http://arxiv.org/abs/0907.0436}

\bibitem{plcJsts07}
Combettes, P.L., Pesquet, J.C.: A {D}ouglas-{R}achford splitting approach to
  nonsmooth convex variational signal recovery.
\newblock IEEE J. Selected Topics Signal Process. \textbf{1}, 564--574 (2007)

\bibitem{plcSiop07}
Combettes, P.L., Pesquet, J.C.: Proximal thresholding algorithm for
  minimization over orthonormal bases.
\newblock SIAM J. Optim. \textbf{18}, 1351--1376 (2007)

\bibitem{plcInvp08}
Combettes, P.L., Pesquet, J.C.: A proximal decomposition method for solving
  convex variational inverse problems.
\newblock Inverse Problems \textbf{24}, 27 (2008).
\newblock Art. 065014

\bibitem{plcSmms05}
Combettes, P.L., Wajs, V.R.: Signal recovery by proximal forward-backward
  splitting.
\newblock Multiscale Model. Simul. \textbf{4}, 1168--1200 (2005)

\bibitem{plcDaub04}
Daubechies, I., Defrise, M., {De Mol}, C.: An iterative thresholding algorithm
  for linear inverse problems with a sparsity constraint.
\newblock Comm. Pure Appl. Math. \textbf{57}, 1413--1457 (2004)

\bibitem{plcDaub07}
Daubechies, I., Teschke, G., Vese, L.: Iteratively solving linear inverse
  problems under general convex constraints.
\newblock Inverse Probl. Imaging \textbf{1}, 29--46 (2007)

\bibitem{plcDemo02}
{De Mol}, C., Defrise, M.: A note on wavelet-based inversion algorithms.
\newblock Contemp. Math. \textbf{313}, 85--96 (2002)

\bibitem{plcDida09}
Didas, S., Setzer, S., Steidl, G.: Combined $\ell_2$ data and gradient fitting
  in conjunction with $\ell_1$ regularization.
\newblock Adv. Comput. Math. \textbf{30}, 79--99 (2009)

\bibitem{plcDoug56}
Douglas, J., Rachford, H.H.: On the numerical solution of heat conduction
  problems in two or three space variables.
\newblock Trans. Amer. Math. Soc. \textbf{82}, 421--439 (1956)

\bibitem{plcDunn76}
Dunn, J.C.: Convexity, monotonicity, and gradient processes in {H}ilbert space.
\newblock J. Math. Anal. Appl. \textbf{53}, 145--158 (1976)

\bibitem{plcDupe09}
Dup\'e, F.X., Fadili, M.J., Starck, J.L.: A proximal iteration for deconvolving
  {P}oisson noisy images using sparse representations.
\newblock IEEE Trans. Image Process. \textbf{18}, 310--321 (2009)

\bibitem{plcDura08}
Durand, S., Fadili, J., Nikolova, M.: Multiplicative noise removal using {L1}
  fidelity on frame coefficients.
\newblock J. Math. Imaging Vision \textbf{36}, 201--226 (2010)

\bibitem{plcDyks83}
Dykstra, R.L.: An algorithm for restricted least squares regression.
\newblock J. Amer. Stat. Assoc. \textbf{78}, 837--842 (1983)

\bibitem{plcEcks94}
Eckstein, J.: Parallel alternating direction multiplier decomposition of convex
  programs.
\newblock J. Optim. Theory Appl. \textbf{80}, 39--62 (1994)

\bibitem{plcEcks92}
Eckstein, J., Bertsekas, D.P.: On the {D}ouglas-{R}achford splitting method and
  the proximal point algorithm for maximal monotone operators.
\newblock Math. Programming \textbf{55}, 293--318 (1992)

\bibitem{plcSvai09}
Eckstein, J., Svaiter, B.F.: General projective splitting methods for sums of
  maximal monotone operators.
\newblock SIAM J. Control Optim. \textbf{48}, 787--811 (2009)

\bibitem{plcEick92}
Eicke, B.: Iteration methods for convexly constrained ill-posed problems in
  {H}ilbert space.
\newblock Numer. Funct. Anal. Optim. \textbf{13}, 413--429 (1992)

\bibitem{plcEsse09}
Esser, E.: Applications of {L}agrangian-based alternating direction methods and
  connections to split {B}regman (2009).
\newblock {ftp://ftp.math.ucla.edu/pub/camreport/cam09-31.pdf}

\bibitem{plcPeyr09}
Fadili, J., Peyr\'e, G.: Total variation projection with first order schemes
  (2009).
\newblock http://hal.archives-ouvertes.fr/hal-00380491

\bibitem{plcFigu09}
Figueiredo, M.A.T., Bioucas-Dias, J.M.: Deconvolution of {P}oissonian images
  using variable splitting and augmented {L}agrangian optimization.
\newblock In: Proc. IEEE Workshop Statist. Signal Process. Cardiff, UK (2009).
\newblock {http://arxiv.org/abs/0904.4872}

\bibitem{plcFigu03}
Figueiredo, M.A.T., Nowak, R.D.: An {EM} algorithm for wavelet-based image
  restoration.
\newblock IEEE Trans. Image Process. \textbf{12}, 906--916 (2003)

\bibitem{plcFigu07}
Figueiredo, M.A.T., Nowak, R.D., Wright, S.J.: Gradient projection for sparse
  reconstruction: Application to compressed sensing and other inverse problems.
\newblock IEEE J. Selected Topics Signal Process. \textbf{1}, 586--597 (2007)

\bibitem{plcForn07}
Fornasier, M.: Domain decomposition methods for linear inverse problems with
  sparsity constraints.
\newblock Inverse Problems \textbf{23}, 2505--2526 (2007)

\bibitem{plcGlow83}
Fortin, M., Glowinski, R. (eds.): Augmented Lagrangian Methods: Applications to
  the Numerical Solution of Boundary-Value Problems.
\newblock North-Holland, Amsterdam (1983)

\bibitem{plcGaba83}
Gabay, D.: Applications of the method of multipliers to variational
  inequalities.
\newblock In: M.~Fortin, R.~Glowinski (eds.) Augmented Lagrangian Methods:
  Applications to the Numerical Solution of Boundary-Value Problems, pp.
  299--331. North-Holland, Amsterdam (1983)

\bibitem{plcGaba76}
Gabay, D., Mercier, B.: A dual algorithm for the solution of nonlinear
  variational problems via finite elements approximations.
\newblock Comput. Math. Appl. \textbf{2}, 17--40 (1976)

\bibitem{plcGlow75}
Glowinski, R., Marrocco, A.: Sur l'approximation, par \'el\'ements finis
  d'ordre un, et la r\'esolution, par p\'enalisation-dualit\'e, d'une classe de
  probl\`emes de {D}irichlet non lin\'eaires.
\newblock RAIRO Anal. Numer. \textbf{2}, 41--76 (1975)

\bibitem{plcGlow89}
Glowinski, R., Tallec, P.L. (eds.): Augmented Lagrangian and Operator-Splitting
  Methods in Nonlinear Mechanics.
\newblock SIAM, Philadelphia (1989)

\bibitem{plcGold85}
Goldburg, M., Marks~II, R.J.: Signal synthesis in the presence of an
  inconsistent set of constraints.
\newblock IEEE Trans. Circuits Syst. \textbf{32}, 647--663 (1985)

\bibitem{plcGold09}
Goldstein, T., Osher, S.: The split {B}regman method for {$L1$}-regularized
  problems.
\newblock SIAM J. Imaging Sci. \textbf{2}, 323--343 (2009)

\bibitem{plcGolu96}
Golub, G.H., {Van Loan}, C.F.: Matrix Computations, 3rd edn.
\newblock Johns Hopkins University Press, Baltimore, MD (1996)

\bibitem{plcGule91}
G{\"u}ler, O.: On the convergence of the proximal point algorithm for convex
  minimization.
\newblock SIAM J. Control Optim. \textbf{20}, 403--419 (1991)

\bibitem{plcGule92}
G{\"u}ler, O.: New proximal point algorithms for convex minimization.
\newblock SIAM J. Optim. \textbf{2}, 649--664 (1992)

\bibitem{plcHale08}
Hale, E.T., Yin, W., Zhang, Y.: Fixed-point continuation for
  $l_1$-minimization: methodology and convergence.
\newblock SIAM J. Optim. \textbf{19}, 1107--1130 (2008)

\bibitem{plcHerm80}
Herman, G.T.: Fundamentals of Computerized Tomography -- Image Reconstruction
  from Projections, 2nd edn.
\newblock Springer-Verlag, London (2009)

\bibitem{plcHiri93}
Hiriart-Urruty, J.B., Lemar\'echal, C.: Fundamentals of Convex Analysis.
\newblock Springer-Verlag, New York (2001)

\bibitem{plcHuan08}
Huang, Y., Ng, M.K., Wen, Y.W.: A fast total variation minimization method for
  image restoration.
\newblock Multiscale Model. Simul. \textbf{7}, 774--795 (2008)

\bibitem{plcJoha06}
Johansson, B., Elfving, T., Kozlovc, V., Censor, Y., Forss\'en, P.E., Granlund,
  G.: The application of an oblique-projected {L}andweber method to a model of
  supervised learning.
\newblock Math. Comput. Modelling \textbf{43}, 892--909 (2006)

\bibitem{plcLopu97}
Kiwiel, K.C., {\L}opuch, B.: Surrogate projection methods for finding fixed
  points of firmly nonexpansive mappings.
\newblock SIAM J. Optim. \textbf{7}, 1084--1102 (1997)

\bibitem{plcLema89}
Lemaire, B.: The proximal algorithm.
\newblock In: J.P. Penot (ed.) New Methods in Optimization and Their Industrial
  Uses, \emph{International Series of Numerical Mathematics}, vol.~87, pp.
  73--87. Birkh\"{a}user, Boston, MA (1989)

\bibitem{plcLema98}
Lemaire, B.: It\'eration et approximation.
\newblock In: \'Equations aux D\'eriv\'ees Partielles et Applications, pp.
  641--653. Gauthiers-Villars, Paris (1998)

\bibitem{plcLent81}
Lent, A., Tuy, H.: An iterative method for the extrapolation of band-limited
  functions.
\newblock J. Math. Anal. Appl. \textbf{83}, 554--565 (1981)

\bibitem{plcLevi66}
Levitin, E.S., Polyak, B.T.: Constrained minimization methods.
\newblock U.S.S.R. Comput. Math. Math. Phys. \textbf{6}, 1--50 (1966)

\bibitem{plcLieu69}
Lieutaud, J.: {A}pproximation d'{O}p\'erateurs par des {M}\'ethodes de
  {D}\'ecomposition.
\newblock Th\`ese, Universit\'e de Paris (1969)

\bibitem{plcLion79}
Lions, P.L., Mercier, B.: Splitting algorithms for the sum of two nonlinear
  operators.
\newblock SIAM J. Numer. Anal. \textbf{16}, 964--979 (1979)

\bibitem{plcLori09}
Loris, I., Bertero, M., {De Mol}, C., Zanella, R., Zanni, L.: Accelerating
  gradient projection methods for $\ell^1$-constrained signal recovery by
  steplength selection rules.
\newblock Appl. Comput. Harm. Anal. \textbf{27}, 247--254 (2009)

\bibitem{plcMart70}
Martinet, B.: R\'egularisation d'in\'equations variationnelles par
  approximations successives.
\newblock Rev. Fran\c{c}aise Informat. Rech. Op\'er. \textbf{4}, 154--158
  (1970)

\bibitem{plcMerc79}
Mercier, B.: Topics in Finite Element Solution of Elliptic Problems.
\newblock No.~63 in Lectures on Mathematics. Tata Institute of Fundamental
  Research, Bombay (1979)

\bibitem{plcMerc80}
Mercier, B.: In\'equations Variationnelles de la M\'ecanique.
\newblock No. 80.01 in Publications Math\'ematiques d'Orsay. Universit\'e de
  Paris-XI, Orsay, France (1980)

\bibitem{plcMor62b}
Moreau, J.J.: Fonctions convexes duales et points proximaux dans un espace
  hilbertien.
\newblock C. R. Acad. Sci. Paris S\'er. A Math. \textbf{255}, 2897--2899 (1962)

\bibitem{plcMore65}
Moreau, J.J.: Proximit\'e et dualit\'e dans un espace hilbertien.
\newblock Bull. Soc. Math. France \textbf{93}, 273--299 (1965)

\bibitem{plcNemi83}
Nemirovsky, A.S., Yudin, D.B.: Problem Complexity and Method Efficiency in
  Optimization.
\newblock Wiley, New York (1983)

\bibitem{plcNest05}
Nesterov, {\relax Yu}.: Smooth minimization of non-smooth functions.
\newblock Math. Program. \textbf{103}, 127--152 (2005)

\bibitem{plcNest07}
Nesterov, {\relax Yu}.: Gradient methods for minimizing composite objective
  function.
\newblock CORE discussion paper 2007076, Universit\'e Catholique de Louvain,
  Center for Operations Research and Econometrics (2007)

\bibitem{plcNest83}
Nesterov, {\relax Yu}.E.: A method of solving a convex programming problem with
  convergence rate $o(1/k^2)$.
\newblock Soviet Math. Dokl. \textbf{27}, 372--376 (1983)

\bibitem{plcNoba95}
Nobakht, R.A., Civanlar, M.R.: Optimal pulse shape design for digital
  communication systems by projections onto convex sets.
\newblock IEEE Trans. Communications \textbf{43}, 2874--2877 (1995)

\bibitem{plcPass79}
Passty, G.B.: Ergodic convergence to a zero of the sum of monotone operators in
  {H}ilbert space.
\newblock J. Math. Anal. Appl. \textbf{72}, 383--390 (1979)

\bibitem{plcPesq96}
Pesquet, J.C., Combettes, P.L.: Wavelet synthesis by alternating projections.
\newblock IEEE Trans. Signal Process. \textbf{44}, 728--732 (1996)

\bibitem{plcPier76}
Pierra, G.: {\'E}clatement de contraintes en parall\`ele pour la minimisation
  d'une forme quadratique.
\newblock Lecture Notes in Comput. Sci. \textbf{41}, 200--218 (1976)

\bibitem{plcPier84}
Pierra, G.: Decomposition through formalization in a product space.
\newblock Math. Programming \textbf{28}, 96--115 (1984)

\bibitem{plcPott93}
Potter, L.C., Arun, K.S.: A dual approach to linear inverse problems with
  convex constraints.
\newblock SIAM J. Control Optim. \textbf{31}, 1080--1092 (1993)

\bibitem{plcPust09}
Pustelnik, N., Chaux, C., Pesquet, J.C.: Parallel proximal algorithm for image
  restoration using hybrid regularization (2009).
\newblock {http://arxiv.org/abs/0911.1536}

\bibitem{plcRock70}
Rockafellar, R.T.: Convex Analysis.
\newblock Princeton University Press, Princeton, NJ (1970)

\bibitem{plcRock76}
Rockafellar, R.T.: Monotone operators and the proximal point algorithm.
\newblock SIAM J. Control Optim. \textbf{14}, 877--898 (1976)

\bibitem{plcRudi92}
Rudin, L.I., Osher, S., Fatemi, E.: Nonlinear total variation based noise
  removal algorithms.
\newblock Physica D \textbf{60}, 259--268 (1992)

\bibitem{plcSetz09}
Setzer, S.: Split {B}regman algorithm, {D}ouglas-{R}achford splitting and frame
  shrinkage.
\newblock Lecture Notes in Comput. Sci. \textbf{5567}, 464--476 (2009)

\bibitem{plcSetz10}
Setzer, S., Steidl, G., Teuber, T.: Deblurring {P}oissonian images by split
  {B}regman techniques.
\newblock J. Vis. Commun. Image Represent. \textbf{21}, 193--199 (2010)

\bibitem{plcSibo70}
Sibony, M.: M\'ethodes it\'eratives pour les \'equations et in\'equations aux
  d\'eriv\'ees partielles non lin\'eaires de type monotone.
\newblock Calcolo \textbf{7}, 65--183 (1970)

\bibitem{plcSpin83}
Spingarn, J.E.: Partial inverse of a monotone operator.
\newblock Appl. Math. Optim. \textbf{10}, 247--265 (1983)

\bibitem{plcStar87}
Stark, H. (ed.): Image Recovery: Theory and Application.
\newblock Academic Press, San Diego, CA (1987)

\bibitem{plcStar98}
Stark, H., Yang, Y.: Vector Space Projections : A Numerical Approach to Signal
  and Image Processing, Neural Nets, and Optics.
\newblock Wiley, New York (1998)

\bibitem{plcStei09}
Steidl, G., Teuber, T.: Removing multiplicative noise by {D}ouglas-{R}achford
  splitting methods.
\newblock J. Math. Imaging Vision \textbf{36}, 168--184 (2010)

\bibitem{plcThom93}
Thompson, A.M., Kay, J.: On some {B}ayesian choices of regularization parameter
  in image restoration.
\newblock Inverse Problems \textbf{9}, 749--761 (1993)

\bibitem{plcTibs96}
Tibshirani, R.: Regression shrinkage and selection via the lasso.
\newblock J. Royal. Statist. Soc. B \textbf{58}, 267--288 (1996)

\bibitem{plcTitt85}
Titterington, D.M.: General structure of regularization procedures in image
  reconstruction.
\newblock Astronom. and Astrophys. \textbf{144}, 381--387 (1985)

\bibitem{plcTrop06}
Tropp, J.A.: Just relax: Convex programming methods for identifying sparse
  signals in noise.
\newblock IEEE Trans. Inform. Theory \textbf{52}, 1030--1051 (2006)

\bibitem{plcTrus84}
Trussell, H.J., Civanlar, M.R.: The feasible solution in signal restoration.
\newblock IEEE Trans. Acoust., Speech, Signal Process. \textbf{32}, 201--212
  (1984)

\bibitem{plcTrus85}
Trussell, H.J., Civanlar, M.R.: The {L}andweber iteration and projection onto
  convex sets.
\newblock IEEE Trans. Acoust., Speech, Signal Process. \textbf{33}, 1632--1634
  (1985)

\bibitem{plcTsen91}
Tseng, P.: Applications of a splitting algorithm to decomposition in convex
  programming and variational inequalities.
\newblock SIAM J. Control Optim. \textbf{29}, 119--138 (1991)

\bibitem{plcTsen08}
Tseng, P.: On accelerated proximal gradient methods for convex-concave
  optimization (2008).
\newblock {http://www.math.washington.edu/$\sim$tseng/papers/apgm.pdf}

\bibitem{plcVarg00}
Varga, R.S.: Matrix Iterative Analysis, 2nd edn.
\newblock Springer-Verlag, New York (2000)

\bibitem{plcVese04}
Vese, L.A., Osher, S.J.: Image denoising and decomposition with total variation
  minimization and oscillatory functions.
\newblock J. Math. Imaging Vision \textbf{20}, 7--18 (2004)

\bibitem{plcVone08}
Vonesh, C., Unser, M.: A fast thresholded {L}andweber algorithm for
  wavelet-regularized multidimensional deconvolution.
\newblock IEEE Trans. Image Process. \textbf{17}, 539--549 (2008)

\bibitem{plcVone09}
Vonesh, C., Unser, M.: A fast multilevel algorithm for wavelet-regularized
  image restoration.
\newblock IEEE Trans. Image Process. \textbf{18}, 509--523 (2009)

\bibitem{plcWeis09}
Weiss, P., Aubert, G., Blanc-F\'eraud, L.: Efficient schemes for total
  variation minimization under constraints in image processing.
\newblock SIAM J. Sci. Comput. \textbf{31}, 2047--2080 (2009)

\bibitem{plcYama98}
Yamada, I., Ogura, N., Yamashita, Y., Sakaniwa, K.: Quadratic optimization of
  fixed points of nonexpansive mappings in {H}ilbert space.
\newblock Numer. Funct. Anal. Optim. \textbf{19}, 165--190 (1998)

\bibitem{plcYoul78}
Youla, D.C.: Generalized image restoration by the method of alternating
  orthogonal projections.
\newblock IEEE Trans. Circuits Syst. \textbf{25}, 694--702 (1978)

\bibitem{plcYoul87}
Youla, D.C.: Mathematical theory of image restoration by the method of convex
  projections.
\newblock In: H.~Stark (ed.) Image Recovery: Theory and Application, pp.
  29--77. Academic Press, San Diego, CA (1987)

\bibitem{plcYoul86}
Youla, D.C., Velasco, V.: Extensions of a result on the synthesis of signals in
  the presence of inconsistent constraints.
\newblock IEEE Trans. Circuits Syst. \textbf{33}, 465--468 (1986)

\bibitem{plcYoul82}
Youla, D.C., Webb, H.: Image restoration by the method of convex projections:
  Part 1 -- theory.
\newblock IEEE Trans. Medical Imaging \textbf{1}, 81--94 (1982)

\bibitem{plcZeidXX}
Zeidler, E.: Nonlinear Functional Analysis and Its Applications, vol. I--V.
\newblock Springer-Verlag, New York (1985--1990)

\bibitem{plcZhan09}
Zhang, X., Burger, M., Bresson, X., Osher, S.: Bregmanized nonlocal
  regularization for deconvolution and sparse reconstruction (2009).
\newblock {ftp://ftp.math.ucla.edu/pub/camreport/cam09-03.pdf}

\bibitem{plcZhuc95}
Zhu, C.Y.: Asymptotic convergence analysis of the forward-backward splitting
  algorithm.
\newblock Math. Oper. Res. \textbf{20}, 449--464 (1995)

\end{thebibliography}

\end{document}